\documentclass{amsart}
\usepackage{amsmath, amssymb}
\usepackage{amsfonts}
\usepackage[arrow,matrix,curve,cmtip,ps]{xy}

\usepackage{amsthm}

\allowdisplaybreaks

\newtheorem{theorem}{Theorem}[section]
\newtheorem{lemma}[theorem]{Lemma}

\newtheorem{corollary}[theorem]{Corollary}

\newtheorem*{theorem*}{Theorem}
\theoremstyle{remark}
\newtheorem{remark}[theorem]{Remark}
\newtheorem{definition}[theorem]{Definition}
\newtheorem{example}[theorem]{Example}

\numberwithin{equation}{section}

\newcommand{\W}{\text{Wojciech}}
\newcommand{\wvec}{\omega}
\newcommand{\wmap}{\omega}
\newcommand{\Z}{\mathbb{Z}}
\newcommand{\N}{\mathbb{N}}

\newcommand{\C}{\mathbb{C}}
\newcommand{\K}{\mathcal{K}}
\newcommand{\B}{\mathcal{B}}
\newcommand{\Q}{\mathcal{Q}}
\newcommand{\M}{\mathcal{M}}
\newcommand{\Hi}{\mathcal{H}}

\newcommand{\im}{\operatorname{im }}
\newcommand{\coker}{\operatorname{coker }}
\newcommand{\Ad}{\operatorname{Ad}}
\newcommand{\ext}{\operatorname{Ext}}
\newcommand{\prim}{\operatorname{Prim}}
\newcommand{\ind}{\operatorname{ind}}
\newcommand{\rank}{\operatorname{rank}}

\begin{document}
\title{Ext classes and Embeddings for $C^*$-algebras of
graphs with sinks}
\author{Mark Tomforde 
}

\address{Department of Mathematics\\ Dartmouth College\\
Hanover\\ NH 03755-3551\\ USA}
\email{mark.tomforde@dartmouth.edu}


\date{\today}

\maketitle

\begin{abstract}
We consider directed graphs $E$ which
have been obtained by adding a sink to a fixed graph $G$. 
We associate an element of $\ext(C^*(G))$ to each such $E$,
and show that the classes of two such graphs are equal in
$\ext(C^*(G))$ if and only if the associated $C^*$-algebra
of one can be embedded as a full corner in the
$C^*$-algebra of the other in a particular way.  If every
loop in $G$ has an exit, then we are able to use this
result to generalize some known classification
theorems for $C^*$-algebras of graphs with sinks.
\end{abstract}

\section{Introduction}

The Cuntz-Krieger algebras $\mathcal{O}_A$ are $C^*$-algebras
generated by a family of partial isometries whose relations
are defined by a finite matrix $A$ with entries in $\{ 0 ,1
\}$ and no zero rows.  In 1982 Watatani \cite{Wat} noted that
one can view $\mathcal{O}_A$ as the $C^*$-algebra of a finite
directed graph $G$ with vertex adjacency
matrix $A$, and the condition that $A$ has no
zero rows implies that $G$ has no sinks.  

In the late 1990's analogues of these
$C^*$-algebras were considered for possibly infinite graphs
which are allowed to contain sinks \cite{KPR, KPRR}.  Since
that time there has been much interest in these graph
algebras.  By allowing graphs which are infinite
and may contain sinks, the class of graph algebras has been
extended to include many $C^*$-algebras besides the
Cuntz-Krieger algebras.  At the same time, graph algebras
remain tractable $C^*$-algebras to study.  Like the
Cuntz-Krieger algebras, their basic structure is
understandable and many of their invariants (such as
$K$-theory or $\ext$) can be readily computed \cite{DT2}. 
Furthermore, it has been found that many results about
Cuntz-Krieger algebras hold for graph algebras with only
minor modifications.

In addition, the graph approach has the advantage
that it provides a convenient tool for visualization.  If
$G$ is a graph and $C^*(G)$ is its associated
$C^*$-algebra, then many facts about $C^*(G)$ can be
translated into properties of $G$ that can be determined by
observation.  Thus $C^*$-algebraic questions may be
translated into (often easier to deal with) graph
questions.  Although a similar thing can be done for
Cuntz-Krieger algebras --- properties of
$\mathcal{O}_A$ can be translated into properties of the
matrix $A$ --- many of these results take a nicer form if one
works with graphs instead of matrices.

Since sinks were specifically excluded from the original
Cuntz-Krieger treatments as well as some of the earlier graph
algebra work, there is now some interest in investigating the
effect of sinks on the structure of graph algebras. 
This interest is further motivated by a desire to
understand the Exel-Laca algebras of \cite{EL}, which may be
thought of as Cuntz-Krieger algebras of infinite matrices. 
It was shown in
\cite{RS} that Exel-Laca algebras can be realized as direct
limits of $C^*$-algebras of finite graphs with sinks. 
Therefore, it is reasonable to believe that results regarding
$C^*$-algebras of graphs with sinks could prove useful in
the study of Exel-Laca algebras.

Some progress in the study of $C^*$-algebras of graphs with
sinks was made in \cite{RTW} where the authors
looked at a fixed graph $G$ and considered \emph{1-sink
extensions} of $G$.  Loosely speaking, a 1-sink extension
$(E,v_0)$ of $G$ is a graph $E$ which is formed by
adding a single sink $v_0$ to $G$.  We say a 1-sink extension
is \emph{essential} if every vertex of $G$ can reach the sink
$v_0$.  Here is an example of a graph $G$ and an essential
1-sink extension $E$ of $G$.
\vspace{.3in}
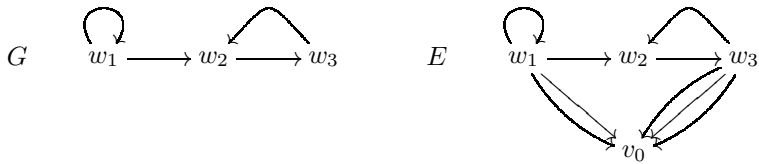
$$ \begin{matrix} G & &
\xymatrix{
w_1 \ar@(ul,ur) \ar[r] & w_2 \ar[r] & w_3 \ar@(ul,ur)[l] \\
}
& & & E & &
\xymatrix{
w_1 \ar@(ul,ur) \ar[r] \ar[dr] \ar@/_/[dr] & w_2 \ar[r]  & w_3
\ar@(ul,ur)[l] \ar[dl] \ar@/_/[dl] \ar@/^/[dl] \\
 & v_0 &  \\
}
\end{matrix}
$$

As in \cite{RTW}, we may associate an
invariant called the \emph{$\W$ vector} to a 1-sink
extension.  This vector is the element $\wvec_E \in
\prod_{G^0} \N$ whose
$w^{\text{th}}$ entry is the number of paths in $E^1
\backslash G^1$ from $w$ to the sink $v_0$.  For instance, the
$\W$ vector in the above example is $\wvec_E = \left(
\begin{smallmatrix} 2 \\ 0 \\ 3
\end{smallmatrix} \right)$.

It was shown in \cite{RTW} that for any 1-sink extension $E$
of $G$ there is an exact sequence 
\begin{equation} \label{sesforext} \xymatrix{ 0 \ar[r] &
I_{v_0} \ar[r]^-i & C^*(E)
\ar[r]^{\pi_E} & C^*(G) \ar[r] & 0. \\}\end{equation} 
Here $I_{v_0}$ denotes the ideal generated by the
projection $p_{v_0}$ corresponding to the sink $v_0$.  If
$E_1$ and $E_2$ are 1-sink extensions, then we say that
$C^*(E_2)$ may be \emph{$C^*(G)$-embedded} into $C^*(E_1)$ if
$C^*(E_2)$ is isomorphic to a full corner of $C^*(E_1)$ via
an isomorphism which commutes with the $\pi_{E_i}$'s.

It was shown in \cite{RTW} that $C^*(G)$-embeddability of
1-sink extensions is determined by the class of the $\W$
vector in $\coker(A_G-I)$, where $A_G$ is the vertex matrix
of $G$.  Specifically, it was shown in
\cite[Theorem 2.3]{RTW} that if $G$ is a graph with no sinks
or sources, $(E_1,v_1)$ and $(E_2,v_2)$ are two essential
1-sink extensions of $G$ whose $\W$ vectors have only a
finite number of nonzero entries, and $\wvec_{E_1}$ and
$\wvec_{E_2}$ are in the same class in $\coker (A_G-I)$, then
there exists a 1-sink extension $F$ of $G$ such that
$C^*(F)$ may be
$C^*(G)$-embedded in both $C^*(E_1)$ and $C^*(E_2)$.  In
addition, a version of this result was proven for
non-essential 1-sink extensions \cite[Proposition 3.3]{RTW}
and a partial converse for both results was obtained in
\cite[Corollary 5.4]{RTW}.  In this paper we show that when
every loop in $G$ has an exit, much stronger results hold.

We shall see in \S \ref{CKiffembed} that if $(E,v_0)$ is
a 1-sink extension of $G$, then (except in degenerate cases)
we will have $I_{v_0} \cong \K$.  Thus we see from
(\ref{sesforext}) that $C^*(E)$ is an extension of $C^*(G)$
by the compact operators.  Hence, $E$ determines an element in
$\ext(C^*(G))$.  In \S \ref{CKiffembed} we prove the
following.  
\begin{theorem*}  Let $G$ be a row-finite graph and
$(E_1,v_1)$ and $(E_2,v_2)$ be 1-sink extensions of $G$.  Then
one of the $C^*(E_i)$'s may be $C^*(G)$-embedded into the
other if and only if $E_1$ and $E_2$ determine the same
element in $\ext(C^*(G))$. 
\end{theorem*}
It was shown in \cite{Tom} that if $G$ is a graph in which
every loop has an exit, then $\ext(C^*(G)) \cong \coker
(A_G-I)$.  Using the isomorphism constructed there we are
able to translate the above result into a statement about the
$\W$ vectors.  Specifically we prove the following.

\begin{theorem*}  Let $G$ be a row-finite graph in which
every loop has an exit.  If $(E_1,v_1)$ and $(E_2,v_2)$
are essential 1-sink extensions of $G$, then one of the
$C^*(E_i)$'s may be $C^*(G)$-embedded into the other if and
only if $[\wvec_{E_1}]=[\wvec_{E_2}]$ in $\coker(A_G-I)$.
\end{theorem*}

Provided that one is willing to allow all the loops in $G$ to
have exits, this result is an improvement over
\cite[Theorem 2.3]{RTW} in the following respects.  First
of all, $G$ is allowed to have sources and there are no
conditions on the $\W$ vectors of $E_1$ and $E_2$.  Second,
we see that the graph $F$ in the statement of \cite[Theorem
2.3]{RTW} can actually be chosen to be either
$E_1$ or $E_2$.  And finally, we see that the equality of the
$\W$ vectors in $\coker (A_G-I)$ is not only sufficient but
necessary.  In \S \ref{nonessembed} we obtain a version of
this theorem for non-essential extensions.

This paper is organized as follows.  We begin in \S
\ref{prelim} with some preliminaries regarding graph
algebras.  We also give precise definitions of 1-sink
extensions, the $\W$ vector, and $C^*(G)$-embeddability.  In
\S \ref{CKiffembed} we show how to associate an element of
$\ext(C^*(G))$ to a (not necessarily essential) 1-sink
extension.  We then prove that if $E_1$ and $E_2$ are 1-sink
extensions of $G$, then one of the $C^*(E_i)$'s may be
$C^*(G)$-embedded into the other if and only if $E_1$ and
$E_2$ determine the same element in $\ext(C^*(G))$.  In \S
\ref{essembed} we recall the definition of the isomorphism
$\wmap : \ext(C^*(G)) \rightarrow \coker (A_G-I)$ from
\cite{Tom} and we prove that for essential 1-sink extensions,
$C^*(G)$-embeddability may be characterized in terms of the
$\W$ vector.  In \S \ref{nonessembed} we discuss non-essential
extensions and again use the isomorphism $\wmap$ to obtain a
characterization of $C^*(G)$-embeddability for arbitrary
1-sink extensions.  We conclude with an example and
some observations

\section{Preliminaries}
\label{prelim}

A (directed) graph $G = (G^0,G^1,r,s)$ consists of
a countable set $G^0$ of vertices, a countable set $G^1$ of
edges, and maps $r,s:G^1 \rightarrow G^0$ which identify the
range and source of each edge.   A vertex
$v \in G^0$ is called a \emph{sink} if
$s^{-1}(v)=\emptyset$ and a \emph{source} if $r^{-1}(v) =
\emptyset$.   We say that a graph is \emph{row-finite} if
each vertex emits only finitely many edges; that is,
$s^{-1}(v)$ is finite for all $v \in G^0$.  All of our
graphs will be assumed to be row-finite.

If $G$ is a row-finite directed graph, a \emph{Cuntz-Krieger
$G$-family} in a $C^*$-algebra is a set of mutually orthogonal
projections $\{ p_v : v \in G^0 \}$ together with a set of
partial isometries $\{ s_e : e \in G^1 \}$ which satisfy the
\emph{Cuntz-Krieger relations} 
$$s_e^* s_e = p_{r(e)} \ \text{for} \  e \in E^1 \ \
\text{and} \
\ p_v =
\sum_{ \{ e : s(e) =v \} } s_e s_e^* \ \text{whenever } v \in G^0
\text{ is not a sink.}$$
\noindent Then $C^*(G)$ is defined to be the $C^*$-algebra
generated by a universal Cuntz-Krieger $G$-family
\cite[Theorem 1.2]{KPR}.

A \emph{path} in a graph $G$ is a finite sequence of edges
$\alpha := \alpha_1 \alpha_2 \ldots \alpha_n$ for which
$r(\alpha_i) = s(\alpha_{i+1})$ for $1 \leq i \leq n-1$, and
we say that such a path has length $|\alpha|=n$.  For $v,w \in
G^0$ we write $v \geq w$ to mean that there exists a path
with source $v$ and range $w$.  For $K,L
\subseteq G^0$ we write $K \geq L$ to mean that for each $v \in
K$ there exists $w \in L$ such that $v \geq w$.  A \emph{loop}
is a path whose range and source are equal, and for a given
loop $x := x_1 x_2 \ldots x_n$ we say that
$x$ is based at $s(x_1) = r(x_n)$.  An
\emph{exit} for a loop $x$ is an edge $e$ for which
$s(e)=s(x_i)$ for some $i$ and $e \neq x_i$.  A graph
is said to satisfy Condition~(L) if every loop in $G$ has an
exit.

We call a loop \emph{simple} if it returns to
its base point exactly once; that is $s(x_1)
\neq r(x_i)$ for $1
\leq i < n$.  A graph is said to satisfy Condition~(K) if no
vertex in the graph is the base of exactly one simple loop;
that is, every vertex is either the base of no loops or the
base of more than one simple loop.  Note that Condition~(K)
implies Condition~(L).

A \emph{$1$-sink extension} of $G$ is a row-finite graph $E$
which contains $G$ as a subgraph and satisfies:
\begin{enumerate}
\item
$H:=E^0\setminus G^0$ is finite, contains no sources,
and contains exactly $1$ sink.
\item
There are no loops in $E$ whose vertices lie in
$H$.
\item
If $e \in E^1 \setminus G^1$, then $r(e) \in H$.
\item
If $w$ is a sink in $G$, then $w$ is a sink in $E$.
\end{enumerate}
When we say $(E, v_0)$ is a $1$-sink extension of $G$, we
mean that
$v_0$ is the sink outside $G^0$.  An edge $e$ with
$r(e)\in H$ and $s(e)\in G^0$ is called a \emph{boundary
edge} and the sources of the boundary edges are called
\emph{boundary vertices}.  We write
$B_E^1$ for the set of boundary edges and $B_E^0$ for the set
of boundary vertices. If
$w\in G^0$ we denote by
$Z(w,v_0)$ the set of paths $\alpha$ from $w$ to $v_0$ which
leave $G$ immediately in the sense that $r(\alpha_1) \notin
G^0$. The \emph{Wojciech vector} of $E$ is the element
$\wvec_E$ of $\prod_{G^0}\N$ given by
$$ \wvec_E(w):=\#Z(w,v_0)\ \text{for $w\in G^0$.}$$

If $(E,v_0)$ is a 1-sink extension of $G$, then there exists
a surjection $\pi_E : C^*(E) \rightarrow C^*(G)$ for which
$$\xymatrix{ 0 \ar[r] & I_{v_0} \ar[r]^-i & C^*(E)
\ar[r]^{\pi_E} & C^*(G) \ar[r] & 0 \\}$$ is a short exact
sequence \cite[Corollary 1.3]{RTW}.  Here $I_{v_0}$ denotes
the ideal in $C^*(E)$ generated by the projection $p_{v_0}$
corresponding to the sink $v_0$.  We say that $(E,v_0)$ is an
\emph{essential} 1-sink extension if $G^0 \geq v_0$.  It
follows from
\cite[Lemma 2.2]{RTW} that $E$ is an essential 1-sink
extension if and only if $I_{v_0}$ is an essential ideal in
$C^*(E)$.  Also note that if there exists an essential 1-sink
extension of $G$, then $G$ cannot have any sinks.

Suppose $(E_1,v_1)$ and $(E_2,v_2)$ are $1$-sink extensions of
$G$. We say that $C^*(E_2)$ is \emph{$C^*(G)$-embeddable into
$C^*(E_1)$} if there is an isomorphism $\phi$ of
$C^*(E_2)=C^*(s_e,p_v)$ onto a full corner in
$C^*(E_1)=C^*(t_f,q_w)$ such that $\phi(p_{v_2})=q_{v_1}$ and
$\pi_E\circ\phi=\pi_F:C^*(F)\to C^*(G)$. We call $\phi$
a \emph{$C^*(G)$-embedding}.  Notice that if $C^*(E_2)$ is
$C^*(G)$-embeddable into $C^*(E_1)$, then
$C^*(E_2)$ is Morita equivalent to $C^*(E_1)$ in a way which
respects the common quotient $C^*(G)$.

If $G$ is a graph, the \emph{vertex matrix} of $G$ is the
$G^0
\times G^0$ matrix
$A_G$ whose entries are given by $A_G(v,w) := \# \{ e \in G^1 :
s(e) = v \text{ and } r(e) = w \}$, and the \emph{edge matrix}
of $G$ is the $G^1 \times G^1$ matrix
$B_G$ whose entries are given by $$B_G(e,f) := \begin{cases}
1 & \text{ if $r(e) = s(f)$.} \\ 0 & \text{
otherwise. } \end{cases}$$  We shall frequently be concerned
with the maps $A_G - I: \prod_{G^0} \Z
\rightarrow \prod_{G^0} \Z$ and $B_G -I : \prod_{G^1} \Z
\rightarrow \prod_{G^1} \Z$ given by left multiplication. 

Throughout we shall let $\Hi$ denote a separable
infinite-dimensional Hilbert space, $\K$ the compact
operators on $\Hi$, $\B$ the bounded operators on $\Hi$, and
$\Q := \B / \K$ the associated Calkin algebra.  We shall also
let $i : \K \rightarrow \B$ denote the inclusion map and $\pi
: \B \rightarrow \Q$ the projection map.  If $A$ is a
$C^*$-algebra, then an \emph{extension} of $A$ (by the
compact operators) is a homomorphism $\tau : A \rightarrow
\Q$.  An extension is said to be \emph{essential} if it is a
monomorphism.

\section{$C^*(G)$-embeddability and CK-equivalence}
\label{CKiffembed}

In order to see how $\ext(C^*(G))$ and
$C^*(G)$-embeddability are related, we will follow the
approach in \cite[\S 3]{Tom} and view $\ext$ as the
CK-equivalence classes of essential extensions.

\begin{definition}
If $\tau_1$ and $\tau_2$ are two (not necessarily essential)
extensions of $A$ by $\K$, then $\tau_1$ and $\tau_2$ are
\emph{CK-equivalent} if there exists either an isometry or
coisometry $W \in \B$ for which
$$\tau_1 = \Ad ( \pi (W) ) \circ \tau_2 \ \text{ and } \ \tau_2
= \Ad ( \pi (W^*) ) \circ \tau_1.$$
\end{definition}
\begin{remark}  In light of \cite[Corollary 5.15]{Tom} we see
that the above definition is equivalent to the one given in
\cite[Definition 3.1]{Tom}.  Also note that CK-equivalence
is not obviously an equivalence relation.  However, for
certain classes of extensions (such as essential extensions)
it has been shown to be an equivalence relation \cite[Remark
3.2]{Tom}.
\end{remark}

Recall that if $E$ is a 1-sink extension of $G$ with
sink $v_0$, then it follows from \cite[Corollary 2.2]{KPR}
that $I_{v_0} \cong \K ( \ell^2 ( E^*(v_0)))$ where $E^*(v_0)
= \{ \alpha \in E^* : r(\alpha) = v_0 \}$.  Thus $I_{v_0}
\cong \K$ when $E^*(v_0)$ contains infinitely many elements,
and $I_{v_0} \cong M_n (\mathbb{C})$ when $E^*(v_0)$ contains
a finite number of elements.  If $G$ has no
sources, then  it is easy to see that $E^*(v_0)$ must have
infinitely many elements, and it was shown in
\cite[Lemma 6.6]{Tom} that if $E$ is an essential 1-sink
extension of $G$, then $E^*(v_0)$ will also have
infinitely many elements.  Consequently, in each of these
cases we will have $I_{v_0} \cong \K$.  Furthermore, one can
see from the proof of \cite[Corollary 2.2]{KPR} that $p_{v_0}$
is a minimal projection in $I_{v_0}$.

\begin{definition}
Let $G$ be a row-finite graph and let $(E,v_0)$ be a
1-sink extension of $G$.  If $I_{v_0}
\cong \K$, (i.e. $E^*(v_0)$ has infinitely many elements), then
choose any isomorphism $i_E : \K \rightarrow I_{v_0}$, and
define the \emph{extension associated to
$E$} to be (the strong equivalence class of) the
Busby invariant $\tau : C^*(G) \rightarrow \Q$ associated to
the short exact sequence 
\begin{equation}
\xymatrix{0 \ar[r] & \K \ar[r]^-{i_E} & C^*(E)
\ar[r]^{\pi_E} & C^*(G) \ar[r] & 0.}
\notag
\end{equation}  

If $I_{v_0} \cong M_n(\mathbb{C})$ for some $n \in \N$ (i.e.
$E^*(v_0)$ has finitely many elements), then the
\emph{extension associated to $E$} is defined to be (the
strong equivalence class of) the zero map $\tau : C^*(G)
\rightarrow \Q$.  That is, $\tau : C^*(G) \rightarrow \Q$
and $\tau (x) = 0$ for all $x \in C^*(G)$.
\end{definition}

Note that the extension associated to $E$ is always a map from
$C^*(G)$ into $\Q$.  Also note that the above definition is
well-defined in the case when $I_{v_0} \cong \K$. That is,
two different choices of $i_E$ will produce extensions with
strongly equivalent Busby invariants (see problem 3E(c) of
\cite{WO} for more details).  Also, since $p_{v_0}$
is a minimal projection, $i_{E}^{-1}(p_{v_0})$ will always
be a rank 1 projection.

Our goal in the remainder of this section is to prove the
following theorem and its corollary.

\begin{theorem}
Let $G$ be a row-finite graph, and let $E_1$ and $E_2$ be
1-sink extensions of $G$. Then the extensions associated
to $E_1$ and $E_2$ are CK-equivalent if and only if one of
the $C^*(E_i)$'s may be $C^*(G)$-embedded into the other.
\label{GequiviffCKequiv}
\end{theorem}

\begin{corollary}
Let $G$ be a row-finite graph, and let $E_1$ and $E_2$ be
essential 1-sink extensions of $G$.  Then the extensions
associated to $E_1$ and $E_2$ are equal in $\ext(C^*(G))$ if
and only if one of the $C^*(E_i)$'s may be $C^*(G)$-embedded
into the other.
\label{Gequiviffwsequiv}
\end{corollary}  

\begin{remark}
Note that we are not assuming that each of the $C^*(E_i)$'s
may be $C^*(G)$-embedded into the other, only that one of
them can.
\end{remark}

\noindent \emph{Proof of Corollary \ref{Gequiviffwsequiv}
.}  Because $E_1$ and $E_2$ are essential it follows from
\cite[Lemma 6.6]{Tom} that $I_{v_1}
\cong I_{v_2} \cong \K$.  Furthermore, \cite[Lemma 3.2]{Tom}
shows that two essential extensions are equal in $\ext$ if and
only if they are CK-equivalent. \hfill $\square$

\begin{lemma}
Let $P$ and $Q$ be rank 1 projections in $\B$.  Then there
exists a unitary $U \in \B$ such that $P =U^*QU$ and $I-U$ has
finite rank.
\label{unitarymovesprojection}
\end{lemma}
\begin{proof} Straightforward.  \end{proof}

\begin{lemma}
Let $G$ be a row-finite graph, and let $(E,v_0)$ be a
1-sink extension of $G$.  If the extension associated to
$E$ is the zero map, then there is an isomorphism $ \Psi :
C^*(E) \rightarrow C^*(G) \oplus I_{v_0}$ which makes the
diagram  
\begin{equation}
  \notag
  \xymatrix{C^*(E) \ar[r]^-{\Psi} \ar[dr]_{\pi_E}&
C^*(G) \oplus I_{v_0} \ar[d]^{p_1} \\
& C^*(G)} 
\end{equation}
commute.  Here $p_1$ is the projection $(a,b) \mapsto a $.
\label{usefulfactfortrivialexts}
\end{lemma}
\begin{proof}  Since the extension associated to $E$ is
zero, one of two things must occur.  If $I_{v_0} \cong
\K$, then $\tau$ is the Busby invariant of $ 0 \rightarrow
I_{v_0} \overset{i}{\rightarrow} C^*(E)
\overset{\pi_E}{\rightarrow} C^*(G) \rightarrow 0$.  If
$I_{v_0} \cong M_n(\C)$, then since $M_n (\mathbb{C} )$
is unital it follows that $\Q ( I_{v_0})  \cong \mathcal{M} (
M_n (\mathbb{C})) / M_n (\mathbb{C}) = 0$ and the Busby
invariant of $ 0 \rightarrow
I_{v_0} \overset{i}{\rightarrow} C^*(E)
\overset{\pi_E}{\rightarrow} C^*(G) \rightarrow 0$ must be
the zero map. In either case, the Busby invariant of the
extension $ 0 \rightarrow
I_{v_0} \overset{i}{\rightarrow} C^*(E)
\overset{\pi_E}{\rightarrow} C^*(G) \rightarrow 0$ is zero.
From \cite[Proposition 3.2.15]{WO} it follows
that $C^*(E) \cong C^*(G) \oplus I_{v_0}$ via the map $\Psi(x)
:= (\pi_E(x), \sigma(x))$, where $\sigma:C^*(E)
\rightarrow I_{v_0}$ denotes the (unique) map for which
$\sigma \circ i$ is the identity. The fact that
$p_1 \circ \Psi = \pi_E$ then follows by checking each on
generators of $C^*(E)$. \end{proof}

\noindent \emph{Proof of Sufficiency in Theorem
\ref{GequiviffCKequiv} .}  Let $E_1$ and $E_2$
are 1-sink extensions of $G$ whose associated extensions are
CK-equivalent.  Also let $v_1$ and $v_2$ denote the sinks of
$E_1$ and $E_2$ and $\tau_1$ and $\tau_2$ be
the extensions associated to $E_1$ and $E_2$.  Consider the
following cases.

\noindent \textsc{Case 1:}  Either $E^*(v_1)$ is finite or
$E^*(v_2)$ is finite.  

Without loss of generality let us assume that
$E^*(v_1)$ is finite and the number of elements in $E^*(v_1)$
is less than or equal to the number of elements in
$E^*(v_2)$.  Then $I_{v_1} \cong M_n(\mathbb{C} )$ for
some finite $n$, and because $I_{v_2} \cong \K (\ell^2
(E^*(v_2)))$ we see that either $I_{v_2} \cong \K$ or
$I_{v_2} \cong M_m ( \mathbb{C} )$ for some finite $m \geq
n$.  In either case we may choose an imbedding $\phi :
I_{v_1} \rightarrow I_{v_2}$ which maps onto a full
corner of $I_{v_2}$.  (Note that since $I_{v_2}$ is simple
we need only choose $\phi$ to map onto a corner, and then that
corner is automatically full.)  Furthermore, since $p_{v_1}$
and $q_{v_2}$ are rank 1 projections, we may choose $\phi$ in
such a way that $\phi (p_{v_1})=q_{v_2}$.  We now define
$\Phi: C^*(G) \oplus I_{v_1} \rightarrow C^*(G) \oplus
I_{v_2}$ by $\Phi((a,b) )=(a,\phi(b))$.  We see that
$\Phi$ maps
$C^*(G) \oplus I_{v_1}$ onto a full corner of $C^*(G) \oplus
I_{v_2}$ and that $\Phi$ makes the diagram
\begin{equation}
  \xymatrix{C^*(G) \oplus I_{v_1} \ar[rr]^{\Phi}
\ar[dr]_{p_1} & &
C^*(G) \oplus I_{v_2} \ar[dl]^{p_1} \\
& C^*(G) & }   
\notag
\end{equation}
commute, where $p_1$ is the projection $(a,b) \mapsto a$.  Now
since
$\tau_1 = 0$ and
$\tau_2$ is CK-equivalent to $\tau_2$, it follows that $\tau_2
= 0$.  Thus Lemma \ref{usefulfactfortrivialexts}, the
existence of $\Phi$, and the above commutative diagram imply
that $C^*(E_1)$ is $C^*(G)$-embeddable into $C^*(E_2)$.

\noindent \textsc{Case 2:}  Both $E^*(v_1)$ and $E^*(v_2)$
are infinite.  

Then $I_{v_1} \cong I_{v_2} \cong \K$.  Let
$\{s_e, p_v \}$ be the canonical Cuntz-Krieger $E_1$-family in
$C^*(E_1)$, and let $\{t_e, q_v \}$ be the canonical
Cuntz-Krieger $E_2$-family in $C^*(E_2)$.  For $k \in \{ 1,
2\}$, choose isomorphisms
$i_{E_k} : \K \rightarrow I_{v_k}$ so that the Busby
invariant $\tau_k$ of
\begin{equation}
\xymatrix{0 \ar[r] & \K \ar[r]^-{i_{E_k}} & C^*(E_k)
\ar[r]^{\pi_{E_k}} & C^*(G) \ar[r] & 0}
\notag
\end{equation}
is an extension associated to $E_k$.  By hypothesis $\tau_1$
and $\tau_2$ are CK-equivalent.  Therefore, after
interchanging the roles of $E_1$ and $E_2$ if necessary, we
may assume that there exists an isometry $W \in \B$ for which
$\tau_1 = \Ad ( \pi (W) )
\circ \tau_2$ and $\tau_2 = \Ad ( \pi (W^*) ) \circ \tau_1$.  

For $k \in \{ 1, 2 \}$, let
$PB_k := \{ (T,a) \in \B \oplus C^*(G) :
\pi(T)=\tau_k(a)
\}$ be the pullback $C^*$-algebra along $\pi$ and $\tau_k$. 
It follows from \cite[Proposition 3.2.11]{WO} that $PB_k \cong
C^*(E_k)$.  Now for $k \in \{1, 2
\}$, let $\sigma_k$ be the unique map which makes the
diagram 
\begin{equation}
  \notag
  \xymatrix{\K \ar[r]^-{i_{E_k}}
\ar[dr]_{i}&
C^*(E_k) \ar[d]^{\sigma_k} \\
& \B} 
\end{equation}
commute.  Then $\sigma_1(p_{v_1})$ and $\sigma_2(q_{v_2})$ are
rank 1 projections in $\B$.  Choose a unit vector $x \in (\ker
W^*)^\perp$.  By Lemma \ref{unitarymovesprojection} there
exists a unitary $U_1 \in \B$ such that $U_1 \sigma_2(q_{v_2})
U_1^*$ is the projection onto $\text{span} \{ x \}$, and for
which $I-U_1$ is compact.  Therefore, by the way in which
$x$ was chosen $WU_1 \sigma_2(q_{v_2}) U_1^* W^*$ is a rank 1
projection.  We may then use Lemma
\ref{unitarymovesprojection} again to produce a unitary $U_2
\in \B$ for which $U_2 (WU_1 \sigma_2(q_{v_2}) U_1^* W^*)U_2^*
= \sigma_1(p_{v_1})$, and $I-U_2$ is compact.

Let $V := U_2WU_1$.  Then $V$ is an isometry, and we may
define a map $\Phi : PB_2 \rightarrow PB_1$ by $\Phi ((T,a)) =
(VTV^*, a)$.  Since $V^*V=I$ it follows that $\Phi$ is a
homomorphism, and since $U_1$ and $U_2$ differ from $I$ by a
compact operator, we see that $\pi(V) = \pi(W)$.  Therefore
$$\pi (VTV^*) = \pi(W) \pi(T) \pi(W^*) = \pi (W)
\tau_2(a) \pi(W^*) = \tau_1(a)$$
so $(VTV^*,a) \in PB_1$, and $\Phi$ does in fact take values
in $PB_1$.

For $k \in \{ 1, 2 \}$, let $p_k : PB_k \rightarrow C^*(G)$
be the projection $p_k((T,a)) = a$.  Then the diagram  
\begin{equation}
  \xymatrix{PB_2 \ar[rr]^{\Phi}
\ar[dr]_{p_2} & &
PB_1 \ar[dl]^{p_1} \\
& C^*(G) & }   
\notag
\end{equation}
commutes and $\Phi ((\sigma_2(q_{v_2}),0)) =
(\sigma_1(p_{v_1}), 0)$.  Also, for $k \in \{ 1, 2 \}$, let
$\Psi_k$ be the standard isomorphism from $C^*(E_k)$ to $PB_k$
given by $\Psi_k (x) = (\sigma_1(x), \pi_{E_k}(x) )$
\cite[Proposition 3.2.11]{WO}.  Then for each $k \in \{1, 2
\}$, the diagram
\begin{equation}
  \notag
  \xymatrix{C^*(E_k) \ar[rr]^{\Psi_k}
\ar[dr]_{\pi_{E_k}} & &
PB_k \ar[dl]^{p_k} \\
& C^*(G) & }   
\end{equation}
commutes and we have that $\Psi_1(p_{v_1}) =
(\sigma_1(p_{v_1}), 0 )$ and
$\Psi_2(q_{v_2}) = (\sigma_2(q_{v_2}), 0 )$.  If we define
$\phi:C^*(E_2) \rightarrow C^*(E_1)$ by $\phi := \Psi_1^{-1}
\circ \Phi \circ \Psi_2$, then the diagram
\begin{equation}
  \xymatrix{C^*(E_2) \ar[rr]^{\phi}
\ar[dr]_{\pi_{E_2}} & &
C^*(E_1) \ar[dl]^{\pi_{E_1}} \\
& C^*(G) & }   
\label{Phicommuteswithproj}
\end{equation}
commutes and $\phi (q_{v_2}) = p_{v_1}$.

We shall now show that $\phi$ embeds $C^*(E_2)$ onto a
full corner of $C^*(E_1)$.  We begin by showing that
$\Phi$ embeds $PB_2$ onto a corner of $PB_1$.  To see that
$\Phi$ is injective, note that since $V$ is an isometry
\begin{align} \| VTV^* \|^2 = & \ \| (VTV^*)(VTV^*)^* \| = \|
VTV^*VT^*V^* \| = \| VTT^*V^* \| \notag \\
= & \ \| (VT)(VT)^* \| = \| VT \|^2 = \| T \|^2. \notag
\end{align}
Therefore $\| VTV^* \| = \| T \|$, and 
$$\| \Phi ((T,a)) \| = \| (VTV^*, a) \| = \text{max} \{ \|
VTV^* \| , \| a \| \} = \text{max} \{ \| T \| , \| a \| \} =
\| (T, a ) \|.$$

Next we shall show that the image of $\Phi$ is a corner in
$PB_1$.  Let $P := VV^*$ be the range projection of $V$. 
We shall define a map $L_P : PB_1 \rightarrow PB_1$ by
$L_P ((T,a)) = (PT, a)$.  To see that $L_P$ actually takes
values in $PB_1$ recall that $U_1$ and $U_2$ differ
from $I$ by a compact operator and therefore $\pi (V) =
\pi (W)$.  We then have that
\begin{align}
\pi(PT) = & \ \pi (VV^*) \pi(T) 
=   \pi (WW^*) \tau_1(a) = \pi(WW^*) \pi(W) \tau_2(a) \pi
(W^*) \notag \\ = & \ \pi(W) \tau_2(a) \pi (W^*)= \
\tau_1(a). \notag
\end{align}
Hence $(PT,a) \in PB_1$.  In a similar way we may define $R_P
: PB_1 \rightarrow PB_1$ by $R_P((T,a))=(TP,a)$.  Since $P$ is
a projection, we see that $L_P$ and $R_P$ are bounded linear
maps.  One can also check that $(L_P,R_P)$ is a double
centralizer and therefore defines an element $\mathcal{P} : =
(L_P, R_P) \in \mathcal{M}(PB_1)$.  Because $P$
is a projection, $\mathcal{P}$ must also be a
projection. Also for any $(T,a) \in PB_1$ we have
that $\mathcal{P} (T,a) = (PT, a)$ and $(T, a) \mathcal{P} =
(TP,a)$.

Now for all $(T,a) \in PB_2$ we have
\begin{align}
\Phi ((T,a)) = & \ (VTV^*, a) = (VV^*VTV^*VV^*, a) \notag \\
= & \ (PVTV^*P, a) = \mathcal{P} (VTV^*, a) \mathcal{P}
= \mathcal{P} \Phi((T,a)) \mathcal{P} \notag
\end{align}
and therefore $\Phi$ maps $PB_2$ into the corner $\mathcal{P}
(PB_1) \mathcal{P}$.  We shall now show that $\Phi$ actually
maps onto this corner.  If $(T,a) \in \mathcal{P} (PB_1)
\mathcal{P}$, then $$\pi(V^*TV) = \pi(W)^* \pi(T) \pi(W) =
\pi(W)^* \tau_1(a) \pi(W) = \tau_2(a)$$ and so $(VTV^*,a) \in
PB_2$.  But then $\Phi((V^*TV,a)) = (VV^*TVV^*,a) = (PTP,a) =
\mathcal{P} (T,a) \mathcal{P} = (T,a)$.  Thus $\Phi$ embeds
$PB_2$ onto the corner $\mathcal{P} (PB_1) \mathcal{P}$.

Because $\Psi_1$ and $\Psi_2$ are isomorphisms, it follows
that $\phi$ embeds $C^*(E_2)$ onto a corner of $C^*(E_1)$.  We
shall now show that this corner must be full.  This will
follow from the commutativity of diagram
\ref{Phicommuteswithproj}.  Let $I$ be any ideal in
$C^*(E_1)$ with the property that $\im \phi \subseteq I$. 
Since $\phi(q_{v_2})=p_{v_1}$ it follows that
$p_{v_1} \in \im \phi \subseteq I$.  Therefore, $I_{v_1}
\subseteq I$.  Furthermore, for any $w \in G^0$ we have by
commutativity that $\pi_{E_1} (p_w - \phi(q_w) ) = 0$. 
Therefore $p_w - \phi (q_w) \in \ker \pi_{E_1} = I_{v_1}$,
and it follows that $p_w - \phi (q_w)
\in I_{v_1} \subseteq I$.  Since
$\phi (q_w) \in \im \phi \subseteq I$, this implies that $p_w
\in I$ for all $w \in G^0$.  Thus $p_w \in I$ for all $w \in
G^0 \cup \{v_1 \}$.  If we let $H := \{ v \in E_1^0 : p_v \in
I \}$, then it follows from \cite[Lemma 4.2]{BPRS} that $H$ is
a saturated hereditary subset of $C^*(E_1)$.  Since we see from
above that $H$ contains $G^0 \cup \{v_1 \}$, and since $E_1$
is a 1-sink extension of $G$, it follows that $H =
E_1^0$.  Therefore $I_H = C^*(E_1)$ and since $I_H \subseteq
I$ it follows that $I = C^*(E_1)$.  Hence $\im \phi$ is a
full corner in $C^*(E_1)$. \hfill $\square$
\\

\noindent \emph{Proof of Necessity in Theorem
\ref{GequiviffCKequiv} .}  Let $E_1$ and $E_2$ be
1-sink extensions of $G$ and suppose that $C^*(E_2)$ is
$C^*(G)$-embeddable into $C^*(E_1)$.  Let $v_1$ and $v_2$
denote the sinks of $E_1$ and $E_2$, respectively.  For $k
\in \{ 1, 2 \}$ let $E^*_k(v_k):= \{ \alpha
\in E_k^* :  r(\alpha) = v_k \}$, and let $\phi : C^*(E_2)
\rightarrow C^*(E_1)$ be a $C^*(G)$-embedding.  Consider the
following cases.

\noindent \textsc{Case 1:}  $E_1^*(v_1)$ is finite.

Then $I_{v_1} \cong
M_n(\C)$ for some finite $n$.  Since $\phi(I_{v_2}) \subseteq
I_{v_1}$, and $I_{v_2} \cong \K (
\ell^2 ( E_2^*(v_2) ) )$, a dimension argument implies that
$E_2^*(v_2)$ must be finite.  Thus if $\tau_1$ and
$\tau_2$ are the extensions associated to $E_1$ and
$E_2$, we have that $\tau_1 = \tau_2 = 0$ so that
$\tau_1$ and $\tau_2$ are CK-equivalent.

\noindent \textsc{Case 2:}  $E_1^*(v_1)$ is infinite.

Then $I_{v_1} \cong \K$.  Choose any
isomorphism $i_{E_1} : \K \rightarrow I_{v_1}$, and let
$\sigma : C^*(E_1) \rightarrow \B$ be the (unique) map which
makes the diagram 
\begin{equation}
  \notag
  \xymatrix{\K \ar[r]^-{i_{E_1}}
\ar[dr]_{i}&
C^*(E_1) \ar[d]^{\sigma} \\
& \B} 
\end{equation}
commute.  If we let $\tau_1$ be the corresponding Busby
invariant, then $\tau_1$ is the extension associated to
$E_1$.

Furthermore, we know that $I_{v_2} \cong \K (H)$, where $H$
is a Hilbert space which is finite-dimensional if $E_2^*(v_2)$
is finite and infinite-dimensional if $E_2^*(v_2)$ is
infinite.  Choose an isomorphism $i_{E_2} : \K (H)
\rightarrow I_{v_2}$.  Then the diagram
\begin{equation}
  \xymatrix{0 \ar[r] & \K(H) \ar[r]^{i_{E_2}} & C^*(E_2)
\ar[r]^{\pi_{E_2}} \ar[d]^{\phi} & C^*(G) \ar[r] \ar@{=}[d]
& 0 \\
0 \ar[r] & \K \ar[r]^{i_{E_1}} \ar@{=}[d] & C^*(E_1)
\ar[r]^{\pi_{E_1}} \ar[d]^\sigma & C^*(G) \ar[r]
\ar[d]^{\tau_1} & 0 \\
0 \ar[r] & \K \ar[r]^i & \B \ar[r]^{\pi} & \Q \ar[r] & 0 }   
\label{bigcomdiagram}
\end{equation}
commutes and has exact rows. 

Let $\{ s_e, p_v \}$ be the canonical Cuntz-Krieger
$E_2$-family in $C^*(E_2)$ and $\{ t_e, q_v \}$ be the
canonical Cuntz-Krieger $E_1$-family in $C^*(E_1 )$. 

We shall now define a bounded linear transformation $U : H
\rightarrow \Hi$.  Since $i_{E_2}^{-1} (p_{v_2})$ is a rank 1
projection, we may write
$i_{E_2}^{-1} (p_{v_2}) = e \otimes e$, where $e$ is a unit
vector in $\im i_{E_2}^{-1} (p_{v_2})$.  Likewise, we may
write
$i_{E_1}^{-1}(q_{v_1}) = f \otimes f$ for some unit vector $f
\in \im i_{E_1}^{-1}(q_{v_1})$.  For convenience of notation
write $\beta : = \sigma \circ \phi \circ i_{E_2}$.  Note
that $\phi(p_{v_2})=q_{v_1}$ implies that $\beta ( e \otimes
e ) = f \otimes f$.  Now for any
$h \in H$ define $$U(h) := \beta (h \otimes e) (f).$$  Then
$U$ is a linear transformation and 
\begin{align}
\langle U(h), U(k) \rangle = \ & \langle \beta ( h \otimes
e) (f) , \beta (k \otimes e)(f) \rangle = \langle \beta ( k
\otimes e)^* \beta ( h \otimes e) (f) , f \rangle \notag \\
= & \ \langle \beta ( \langle h,k \rangle (e \otimes e)) (f) ,
f \rangle = \langle h,k \rangle \langle \beta (e \otimes e)
(f) , f \rangle \notag \\
= & \ \langle h,k \rangle \langle (f \otimes f)
(f) , f \rangle = \langle h,k \rangle \langle f , f \rangle =
\langle h,k
\rangle. \notag
\end{align}
Therefore $U$ is an isometry.  

Now since $\phi$ embeds $C^*(E_2)$ onto a full corner of
$C^*(E_1)$, it follows that there exists a projection $p \in
\mathcal{M}(C^*(E_1))$ such that $\im \phi = pC^*(E_1)p$. 
Because $\sigma$ is a nondegenerate representation (since
$\sigma(I_{v_1}) = \K$), it extends to a representation
$\overline{\sigma} : \mathcal{M} (C^*(E_1)) \rightarrow \B$ by
\cite[Corollary 2.51]{RW}.  Let $P :=
\overline{\sigma}(p)$.  We shall show that $\im P
\subseteq \im U$.  Let $g \in \im P$.  Also let $f$ be as
before.  Then $g \otimes f \in \K$ and 
$$\sigma ( p i_{E_1} (g \otimes f)p) =  \overline{\sigma} (p)
\sigma (i_{E_1} (g \otimes f)) \overline{\sigma} (p) = P
(g \otimes f) P = \sigma (i_{E_1} (g \otimes f)).$$  Now since
$p i_{E_1}(g \otimes f) p \in p C^*(E_1) p = \im \phi$,
there exists $a \in C^*(E_2)$ such that $\phi(a) = p
i_{E_1} (g \otimes f ) p$.  In addition, since $\pi_{E_1} :
C^*(E_1) \rightarrow C^*(G)$ is surjective, it extends to a
homomorphism $\overline{\pi}_{E_1} : \mathcal{M} (C^*(E_1))
\rightarrow \mathcal{M}(C^*(G))$ by \cite[Corollary
2.51]{RW}.  By commutativity and exactness we then have that
$$\pi_{E_2}(a)=\pi_{E_1}(\phi(a))=\pi_{E_1}(p i_{E_1}(g
\otimes f) p) = \overline{\pi}_{E_1}(p) \pi_{E_1}(i_{E_1}(g
\otimes f)) \overline{\pi}_{E_1}(p) = 0.$$
Thus $a \in \im i_{E_2}$ by exactness, and we have that
$a=i_{E_2}(T)$ for some $T \in \K (H)$.  Let $h := T(e)$. 
Then
\begin{align}
U(T(e)) = & \ \beta(T(e) \otimes e)(f) = \beta (T \circ (e
\otimes e))(f) = \beta(T) \beta(e \otimes e)
(f) \notag \\
= & \ \beta(T) (f \otimes f) (f) = \sigma (p i_{E_1} (g
\otimes f) p) (f) =
\sigma(i_{E_1}(g \otimes f)) (f) \notag \\
= & \ (g \otimes
f)(f) = \langle f, f \rangle g = g. \notag
\end{align}
Thus $g \in \im U$ and $\im P \subseteq \im U$.

Now if $H$ is a finite-dimensional
space, it follows that $\im U$ is finite-dimensional.  Since
$\im P \subseteq \im U$, this implies that
$P$ has finite rank and hence $\pi (P) = 0$.  Now if $x \in
C^*(G)$, then since $\pi_{E_2}$ is surjective there exists an
element $a
\in C^*(E_2)$ for which $\pi_{E_2}(a) = x$.  Since
$\pi_{E_1} (\phi (a)) = \pi_{E_2}(a) = x$, it follows that
$\tau_1(x) = \pi(\sigma(\phi(a)))$.  But since $\phi(a)
\in \im \phi = p C^*(E_1)p$ we have that $ \phi(a) = p
\phi(a)$ and thus $\tau_1(x) = \pi(\overline{\sigma}(p)
\sigma(\phi(p)) = 0$.  Since $x$ was arbitrary this implies
that $\tau_1 = 0$.  Furthermore, since $H$ is
finite-dimensional, the extension associated to $E_2$ is
$\tau_2=0$.  Thus $\tau_1$ and $\tau_2$ are CK-equivalent.

Therefore, all that remains is to consider the case when $H$
is infinite-dimensional.  In this case $H = \Hi$ and $\K(H) =
\K$.  Furthermore, if $S$ is any
element of $\K$, then for all $h \in \Hi$ we have that 
$$ (\beta(S) \circ U)(h) = \beta(S)(\beta (h \otimes
e)(f)) = \beta(Sh \otimes e)(f) = U(Sh). $$
Since $U$ is an isometry this implies that $U^*\beta(S)U=S$
for all $S \in \K$.  Therefore, $\Ad(U^*)\circ \beta$ is the
inclusion map $i : \K \rightarrow \B$.  Since $\Ad(U^*) \circ
\beta = \Ad(U^*) \circ \sigma \circ \phi \circ i_{E_2}$, this
implies that $\Ad(U^*) \circ \sigma \circ \phi$ is the
unique map which makes the following diagram commute.
\begin{equation}
  \notag
  \xymatrix{\K \ar[r]^-{i_{E_2}}
\ar[dr]_{i}&
C^*(E_2) \ar[d]^{\Ad(U^*) \circ \sigma \circ \phi} \\
& \B} 
\end{equation}
Therefore, if $\tau_2$ is (the Busby invariant of) the
extension associated to $C^*(E_2)$, then by definition
$\tau_2$ is equal to the following.  For any $x \in C^*(G)$
choose an $a \in C^*(E_2)$ for which $\pi_{E_2}(a)=x$.  Then
$\tau_2(x) := \pi(\Ad(U^*) \circ \sigma \circ \phi(a))$. 
Using the commutativity of diagram \ref{bigcomdiagram}, this
implies that
\begin{align}
\tau_2(x) = & \ \Ad(\pi(U^*)) \circ \pi(\sigma(\phi(a)))
= \Ad(\pi(U^*)) \circ \tau_1(\pi_{E_1}(\phi(a))) \notag
\\ 
= & \ \Ad(\pi(U^*)) \circ \tau_1(\pi_{E_2}(a)) =
\Ad(\pi(U^*)) \circ \tau_1(x). \notag
\end{align}
So for all $x \in C^*(G)$ we have that 
\begin{equation}
\tau_2(x) = \pi(U^*) \tau_1(x) \pi(U).
\label{firsthalfCK}
\end{equation}Now if $a$ is any element of $C^*(E_2)$, then
$\phi(a) \in pC^*(E_1)p$.  Thus $\phi(a) = p \phi(a)$ and
$$\sigma ( \phi(a)) = \sigma (p \phi(a)) =
\overline{\sigma}(p) \sigma(\phi(a)) = P \sigma(\phi(a)).$$
Hence $\im \sigma(\phi(a)) \subseteq \im P \subseteq \im U$,
and we have that
$$UU^* \sigma \phi (a) = \sigma \phi (a)
\hspace{.2in} \text{ for all $a \in C^*(E_2).$}$$
Furthermore, for any $x \in C^*(G)$, we may choose an $a \in
C^*(E_2)$ for which $\pi_{E_2}(a)=x$, and using the
commutativity of diagram \ref{bigcomdiagram} we then have
that 
\begin{align}
UU^*\sigma \phi (a) = & \ \sigma \phi(a) \notag \\
\pi(UU^*) \pi \sigma \phi (a) = & \ \pi \sigma \phi(a)
\notag \\ 
\pi(UU^*) \tau_1 \pi_{E_1} \phi (a) = & \ \tau_1 \pi_{E_1}
\phi(a) \notag \\ 
\pi(UU^*) \tau_1 \pi_{E_2} (a) = & \ \tau_1 \pi_{E_2} (a) 
\notag \\ 
\pi(UU^*) \tau_1 (x) = & \ \tau_1 (x). \notag
\end{align}
In addition, this implies that for any $x \in C^*(G)$ we
have that $\pi(UU^*)\tau_1(x^*) = \tau_1(x^*)$, and taking
adjoints this gives that
$$\tau_1(x) \pi(UU^*) = \tau_1(x) \hspace{.2in} \text{ for all
$x \in C^*(G)$. }$$
Thus for all $x \in C^*(G)$ we have
$$ \tau_1(x) = \pi(UU^*) \tau_1(x) \pi(UU^*) = \pi(U) \big(
\pi(U^*) \tau_1(x) \pi (U) \big) \pi(U^*)
= \pi (U) \tau_2(x) \pi (U^*). $$  This, combined with
Eq.(\ref{firsthalfCK}), implies that
$\tau_1 = \Ad(\pi(U)) \circ \tau_2$ and $\tau_2 =
\Ad(\pi(U^*)) \circ \tau_1$.  Since $U$ is an
isometry, $\tau_1$ and
$\tau_2$ are CK-equivalent. \hfill $\square$

\section{$C^*(G)$-embeddability for Essential 1-sink
extensions}
\label{essembed}

In the previous section it was shown that if $E_1$ and
$E_2$ are two 1-sink extensions of $G$, then one of
the $C^*(E_i)$'s can be $C^*(G)$-embedded into
the other if and only if their associated extensions are
CK-equivalent.  While this gives a characterization of
$C^*(G)$-embeddability, it is somewhat unsatisfying due to the
fact that CK-equivalence of the Busby invariants is not an
easily checkable condition.  We shall use the $\W$ map
defined in \cite{Tom} to translate this result into a
statement about the $\W$ vectors of $E_1$ and $E_2$.  We
shall do this for essential 1-sink extensions in this
section, and in the next section we shall consider
non-essential 1-sink extensions.

We begin by recalling the definition of the $\W$ map.  If $E
\in Q$ is a projection, and $X$ is an element of $\Q$ such
that $EXE$ is invertible in $E \Q E$, then we denote by
$\ind_E(X)$ the Fredholm index of $E'X'E'$ in $\im E'$,
where $E'$ is any projection in $\B$ for which
$\pi (E') = E$ and $X$ is any element of $\B$ such that
$\pi(X')=X$.

Let $G$ be a row-finite graph with no sinks which satisfies
Condition~(L), and let $\{ s_e,p_v \}$ be the generating
Cuntz-Krieger $G$-family for $C^*(G)$.  If $\tau : C^*(G)
\rightarrow \Q$ is an essential extension of $C^*(G)$, define
$E_e : = \tau(s_es_e^*)$ for all $e \in G^1$.  If $t : C^*(G)
\rightarrow \Q$ is another essential extension of $C^*(G)$
with the property that $t(s_es_e^*) = E_e$ for all $e \in
G^1$, then we define a vector $d_{\tau, t} \in \prod_{G^1}
\Z$ by $$d_{\tau,t} := -\ind_{E_e} \tau(s_e) t(s_e^*).$$

\noindent We then define the \emph{Cuntz-Krieger}
map $d : \ext (C^*(G)) \rightarrow \coker (B_G-I)$ by $$d (
\tau) := [ d_{\tau,t} ],$$ where $t$ is any degenerate
essential extension of $C^*(G)$ with the property that
$t(s_es_e^*) = \tau(s_es_e^*)$ for all $e \in G^1$.  

Furthermore, we define the \emph{source matrix} of $G$ to be
the $G^0 \times G^1$ matrix $S_G$ defined by $$S_G(v,e) = \begin{cases} 1 &
\text{if $s(e)=v$} \\ 0 & \text{otherwise.} \end{cases}$$  It
follows from \cite[Lemma 6.2]{Tom} that $S_G : \prod_{G^1} \Z
\rightarrow \prod_{G^0} \Z$ induces an
isomorphism $\overline{S_G} : \coker (B_G-I) \rightarrow \coker
(A_G-I)$, and we define the \emph{$\W$ map} $\wmap :
\ext (C^*(G)) \rightarrow \coker (A_G-I)$ by $$\wmap ( \tau )
= \overline{S_G} \circ d.$$  It was shown in \cite[Theorem
6.16]{Tom} that both the Cuntz-Krieger map and the $\W$ map
are isomorphisms.

\begin{theorem}
Let $G$ be a row-finite graph which satisfies Condition~(L). 
Also let $E_1$ and $E_2$ be essential 1-sink extensions of
$G$.  Then one of the $C^*(E_i)$'s may be $C^*(G)$-embedded
into the other if and only if $$ [ \wvec_{E_1} ] = [
\wvec_{E_2} ] \text{ in $\coker (A_G-I)$,} $$ 
where $\wvec_{E_i}$ is the $\W$ vector of $E_i$ and
$A_G-I :
\prod_{G^0} \Z \rightarrow \prod_{G^0} \Z$.
\label{charactofGequivforessential}
\end{theorem}  
\begin{proof}  Let $\tau_1$ and $\tau_2$ be the extensions
associated to $E_1$ and $E_2$, respectively.  It follows from
Corollary \ref{Gequiviffwsequiv} that $\tau_1$ and
$\tau_2$ are in the same equivalence class in $\ext(C^*(G))$
if and only if one of the $C^*(E_i)$'s may be
$C^*(G)$-embedded into the other.  Since $E_1$ and $E_2$ are
essential 1-sink extensions of $G$, the graph $G$
contains no sinks.  By
\cite[Theorem 6.16]{Tom} the $\W$ map $\wmap : \ext (C^*(G))
\rightarrow
\coker (A_G-I)$ is an isomorphism, and by
\cite[Proposition 6.11]{Tom} the value of the $\W$ map on
$\tau_i$ is the class $[\wvec_{E_i} ]$ in $\coker (A_G-I)$.
\end{proof}

\section{$C^*(G)$-embeddability for Non-essential 1-sink
extensions}
\label{nonessembed}

Recall from \cite[\S6]{BPRS} that a
\emph{maximal tail} in a graph
$E$ is a nonempty subset of $E^0$ which is cofinal under $\geq$,
is backwards hereditary ($v\geq w$ and $w\in \gamma$ imply $v\in
\gamma$), and contains no sinks (for each $w\in \gamma$, there
exists $e\in E^1$ with $s(e)=w$ and $r(e)\in \gamma$).  The
set of all maximal tails of $G$ is denoted by $\chi_G$.

Also recall from \cite[\S3]{RTW} that if $(E,v_0)$ is a 1-sink
extension of $G$, then the \emph{closure}
of $v_0$ is the set $$ \overline{v_0}:=\bigcup\{\gamma:\gamma
\text{ is a maximal tail in $G$ and }\gamma\geq
v_0\}.$$ Notice first that the extension is essential if and
only if $\overline{v_0}=G^0$.  Also notice that the closure is
a subset of $G^0$ rather than
$E^0$.  It has been defined in this way so that one may
compare the closures in different extensions.

As in \cite{RTW}, we
mention briefly how this notion of closure is related to  the
closure of sets in $\prim C^*(E)$, as described in
\cite[\S6]{BPRS}. For each
sink $v$, let $\lambda_v:=\{w\in E^0:w\geq v\}$, and let
$$ \Lambda_E:=\chi_E \cup\{\lambda_v:v
\text{ is a sink in $E$}\}.$$  The set $\Lambda_E$ has a
topology in which the closure of a subset $S$ is
$\{\lambda:\lambda\geq\bigcup_{\chi\in S}\chi\}$, and it is
proved in \cite[Corollary 6.5]{BPRS} that when $E$ satisfies
Condition~(K) of \cite{KPRR}, $\lambda \mapsto
I(E^0\setminus\lambda)$ is a homeomorphism of $\Lambda_E$ onto
$\prim C^*(E)$. If $(E,v_0)$ is a 1-sink extension of $G$, then
the only loops in $E$ are those in $G$, so $E$ satisfies
Condition~(K) whenever
$G$ does. A subset of $G^0$ is a maximal tail in $E$ if and only
if it is a maximal tail in $G$, and because every sink in $G$ is
a sink in $E$, we deduce that
$\Lambda_E=\Lambda_G\cup \{\lambda_{v_0}\}$.

We now return to the problem of proving an analogue of
Theorem~\ref{charactofGequivforessential} for non-essential
extensions.

\begin{lemma}  Let $G$ be a graph which satisfies
Condition (K), and let $(E_1,v_1)$ and $(E_2,v_2)$ be 1-sink
extensions of $G$.  If $C^*(E_2)$ is $C^*(G)$-embeddable into
$C^*(E_1)$, then $\overline{v_1} = \overline{v_2}$.
\label{primconverse}
\end{lemma}

\begin{proof}  Let $\phi : C^*(E_2) \rightarrow C^*(E_1)$
be a $C^*(G)$-embedding.  Also let $p \in \M (C^*(E_1))$ be
the projection which determines the full corner $\im \phi$. 
Now for $i \in \{ 1, 2 \}$ we have that $\Lambda_{E_i} =
\Lambda_G \cup \{ \lambda_{v_i} \}$ is homeomorphic to
$\prim C^*(E_i)$ via the map $\lambda \mapsto I_{H_\lambda}$,
where $H_\lambda := E_i^0 \backslash \lambda$
by \cite[Corollary 6.5]{BPRS}.  Furthermore, since $\phi$
embeds $C^*(E_2)$ onto a full corner of $C^*(E_1)$ it follows
that $C^*(E_2)$ is Morita equivalent to $C^*(E_1)$ and the
Rieffel correspondence is a homeomorphism between $\prim
C^*(E_2)$ and $\prim C^*(E_1)$, which in this case is given by
$I \mapsto \phi^{-1}(pIp)$ \cite[Proposition 3.24]{RW}.
Composing the homeomorphisms which we have described, we
obtain a homeomorphism from $h : \Lambda_{E_2}
\rightarrow \Lambda_{E_1}$, where
$h(\lambda)$ is the unique element of $\Lambda_{E_1}$ for
which $\phi(I_{H_\lambda}) = pI_{H_{h(\lambda)}}p$.

We shall now show that this homeomorphism $h$ is equal to
the map $h$ described in \cite[Lemma 3.2]{RTW}; that is $h$
restricts to the identity on $\Lambda_G$.  Let $\lambda \in
\Lambda_G \subseteq \Lambda_{E_2}$.  We begin by showing that
$h(\lambda) \in \Lambda_G$.  Let $\{ s_e,p_v
\}$ be the canonical Cuntz-Krieger $E_2$-family, and let $\{
t_f, q_w
\}$ be the canonical Cuntz-Krieger $E_1$-family.  Since
$\lambda \in \Lambda_G$ it follows that $v_2 \notin \lambda$. 
Therefore $v_2 \in H_\lambda$ and $p_{v_2} \in
I_{H_\lambda}$.  Consequently,
$\phi(p_{v_2}) \in \phi (I_{H_\lambda})$, and since
$\phi(p_{v_2}) = q_{v_1}$ and $pI_{H_{h(\lambda)}}p =
\phi(I_{H_\lambda})$ it follows that $q_{v_1} \in
pI_{H_{h(\lambda)}}p \subseteq I_{H_{h(\lambda)}}$.  Thus $v_1
\in H_{h(\lambda)}$ and $v_1 \notin h(\lambda)$.  It follows
that $h(\lambda) \neq \lambda_{v_1}$, and hence $h(\lambda)
\in \Lambda_G$.

We shall now proceed to show that $h(\lambda)=\lambda$.  Since
$h(\lambda) \in \Lambda_G$ it follows that
$H_{v_1} \subseteq H_{h(\lambda)}$.  Thus $\ker
\pi_{E_1} = I_{H_{v_1}} \subseteq I_{H_{h(\lambda)}}$.  Now
let $w \in \lambda$.  If we let $\{ u_g, r_x \}$ be the
canonical Cuntz-Krieger $G$-family, then  since $w \in G^0$
we have that $\pi_{E_2} (p_w) = r_w$.  It then follows that
$$\pi_{E_1} (\phi (p_w) - q_w) = \pi_{E_1} (\phi (p_w) ) -
\pi_{E_1} (q_w) = \pi_{E_2} (p_w) - \pi_{E_1} (q_w) = r_w -
r_w = 0.$$  Thus $\phi (p_w) - q_w \in \ker \pi_{E_1}
\subseteq I_{H_{h(\lambda)}}$.  We shall now show that $q_w
\notin I_{H_{h(\lambda)}}$.  To do this we suppose that
$q_w \in I_{H_{h(\lambda)}}$ and arrive at a contradiction. 
If $q_w \in I_{H_{h(\lambda)}}$, then we would have that
$\phi(p_w) \in I_{H_{h(\lambda)}}$.  Thus $p \phi (p_w) p \in
p I_{H_{h(\lambda)}} p$ and $p \phi (p_w) p \in \phi (
I_{H_{\lambda}} )$.  Now
$\phi(p_w) \in \phi(C^*(E_2))$ and
$\phi ( C^*(E_2)) = pC^*(E_1)p$.  Hence $p \phi(p_w) p =
\phi(p_w)$ and we have that $\phi (p_w) \in \phi
(I_{H_\lambda})$.  Since $\phi$ is injective this implies that
$q_w \in I_{H_\lambda}$ and $w \in H_\lambda$ and $w \notin
\lambda$ which is a contradiction.  Therefore we must have
that $q_w \notin I_{H_{h(\lambda)}}$ and $w \notin
H_{h(\lambda)}$ and $w \in h(\lambda)$.  Hence $\lambda
\subseteq h(\lambda)$.  

To show inclusion in the other direction let $w \in
h(\lambda)$.  Then $w \in H_{h(\lambda)}$ and $q_w \notin
I_{H_{h(\lambda)}}$.  As above, it is the case that
$\phi(p_w)-q_w \in I_{H_{h(\lambda)}}$.  Therefore, $\phi(p_w)
\notin I_{H_{h(\lambda)}}$ and since
$pI_{H_{h(\lambda)}}p \subseteq I_{H_{h(\lambda)}}$ it follows
that $\phi(p_w) \notin pI_{H_{h(\lambda)}}p$ or $\phi(p_w)
\notin \phi(I_{H_\lambda})$.  Thus $p_w \notin I_{H_\lambda}$
and $w \notin H_\lambda$ and $w \in \lambda$.  Hence
$h(\lambda) \subseteq \lambda$.

Thus $\lambda = h(\lambda)$ for any
$\lambda \in \Lambda_G$, and the map $h : \Lambda_{E_2}
\rightarrow \Lambda_{E_1}$ restricts to the identity on
$\Lambda_G$.  Since this map is a bijection it must
therefore take $\lambda_{v_2}$ to $\lambda_{v_1}$.  Therefore
$h$ is precisely the map described in \cite[Lemma 3.2]{RTW},
and it follows from \cite[Lemma 3.2]{RTW} that $\overline{v_1}
= \overline{v_2}$. \end{proof}

\begin{definition}
Let $G$ be a row-finite graph which satisfies Condition~(K). 
If $(E,v_0)$ is a 1-sink extension of $G$ we define $$H_E :=
G^0 \backslash \overline{v_0}.$$  We call $H_E$ the
\emph{inessential part of $E$}.
\end{definition}

\begin{lemma}
Let $G$ be a row-finite graph which satisfies
Condition~(K) and let $(E,v_0)$ be a 1-sink extension of $G$.
Then $H_E$ is a saturated hereditary subset of $G^0$.
\label{HEsaturatedhereditary}
\end{lemma}

\begin{proof} Let $v \in H_E$ and $e \in G^1$ with $s(e) =
v$.  If
$r(e) \notin H_E$, then $r(e) \in \overline{v_0}$ and hence
$r(e) \in
\gamma$ for some $\gamma \in \chi_G$ with the property
that $\gamma \geq v_0$.  Since maximal tails are
backwards hereditary this implies that $v = s(e) \in
\gamma$.  Hence $v \in
\overline{v_0}$ and $v
\notin H_E$ which is a contradiction.  Thus we must have $r(e)
\in H_E$ and $H_E$ is hereditary.

Suppose that $v \notin H_E$.  Then $v \in \overline{v_0}$ and
$v \in \gamma$ for some $\gamma \in \chi_G$ with the property
that $\gamma \geq v_0$.  Since maximal tails contain no sinks
there exists an edge $e \in G^1$ with $s(e)=v$ and $r(e) \in
\gamma$.  Thus $r(e) \in \overline{v_0}$ and $r(e) \notin
H_E$.  Hence $H_E$ is saturated. \end{proof}

\begin{remark}
Recall that if $A$ is a $C^*$-algebra, then there is a
lattice structure on the set of ideals of $A$ given by $I
\wedge J := I
\cap J$ and $I \vee J :=$ the smallest ideal containing $I
\cup J$.  Furthermore, if $G$ is a graph then the set of
saturated hereditary subsets of $G^0$ also has a lattice
structure given by $H_1 \wedge H_2 := H_1 \cap H_2$ and $H_1
\vee H_2 :=$ the smallest saturated hereditary subset
containing $H_1 \cup H_2$.  If $G$ is a row-finite graph
satisfying Condition~(K), then it is shown in
\cite[Theorem 4.1]{BPRS} that the map $H \mapsto I_H$, where
$I_H$ is the ideal in $C^*(G)$ generated by
$\{p_v : v \in H \}$, is a lattice isomorphism from the
lattice of saturated hereditary subsets of $G^0$ onto the
lattice of ideals of $C^*(G)$.  We shall make use of this
isomorphism in the following lemmas in order to calculate
$\ker \tau$ for an extension $\tau : C^*(G) \rightarrow \Q$.
\end{remark}

\begin{lemma}
Let $0 \rightarrow \K \overset{i_E}{\rightarrow} E
\overset{\pi_E}{\rightarrow} A \rightarrow 0$ be a short
exact sequence, and let $\sigma$ and $\tau$ be the unique maps
which make the diagram
\begin{equation}
  \notag
  \xymatrix{0 \ar[r] & \K \ar@{=}[d] \ar[r]^{i_E} &
E \ar[d]^{\sigma} \ar[r]^{\pi_E} & A \ar[r] \ar[d]^{\tau} & 0
\\  0 \ar[r] & \K \ar[r]^i &
\B \ar[r]^{\pi} & \Q \ar[r] & 0}
\end{equation}
commute.  Then $\ker (\pi \circ \sigma) = i_E(\K) \vee \ker
\sigma$ and $\ker \tau = \pi_E ( i_E (\K) \vee \ker \sigma )$.
\label{kerneloftauisthewedgeofideals}
\end{lemma}

\begin{proof}  Since $\ker (\pi \circ \sigma)$ is an ideal
which contains $i_E(K)$ and $\ker \sigma$, it follows that
$i_E(\K)
\vee \ker \sigma \subseteq \ker (\pi \circ \sigma)$.

Conversely, if $x \in \ker (\pi \circ \sigma)$ then
$\pi(\sigma(x))=0$ and $\sigma(x) \in \K = \sigma(i_E(\K))$. 
Thus $\sigma (x) = \sigma(a)$ for some $a \in i_E(\K)$.  Hence
$x-a \in \ker \sigma$ and $x \in i_E(\K) \vee \ker \sigma$. 
Thus $\ker (\pi \circ \sigma) = i_E(\K) \vee \ker \sigma$.

In addition, the commutativity of the above diagram implies
that $\pi^{-1} (\ker \tau) = \ker (\pi \circ \tau)$.  Since
$\pi_E$ is surjective it follows that $\ker \tau = \pi_E (
\ker ( \pi \circ \sigma) )$ and from the previous paragraph
$\ker \tau = \pi_E (i_E(\K) \vee \ker \sigma)$. \end{proof}

For Lemmas \ref{kernelofsigmaisvdontgotosink} and
\ref{wedgeequalsunion} fix a row-finite graph $G$ which
satisfies Condition~(K).  Also let $(E,v_0)$ be a fixed
1-sink extension of $G$ which has the property that
$E^*(v_0) := \{ \alpha \in E^* : r(\alpha) = v_0 \}$ contains
infinitely many elements.  Then $I_{v_0} \cong \K$, and we
may choose an isomorphism $i_E : \K \rightarrow I_{v_0}$ and
let $\sigma$ and $\tau$ be the (unique) maps which make the
diagram
\begin{equation}
  \notag
  \xymatrix{0 \ar[r] & \K \ar@{=}[d] \ar[r]^{i_E} &
C^*(E) \ar[d]^{\sigma} \ar[r]^{\pi_E} & C^*(G) \ar[r]
\ar[d]^{\tau} & 0
\\  0 \ar[r] & \K \ar[r]^i &
\B \ar[r]^{\pi} & \Q \ar[r] & 0 \\}
\end{equation}
commute.  In particular, note that $\tau$ is the extension
associated to $E$.

\begin{lemma}
If $\sigma$ is as above, then $\ker \sigma = I_{H'}$ where
$H' := \{ v \in E^0 : v \ngeq v_0 \}$.
\label{kernelofsigmaisvdontgotosink}
\end{lemma}
\begin{proof}  Since $G$ satisfies Condition~(K) and $E$ is a
1-sink extension of $G$, it follows that $E$ also satisfies
Condition~(K).  Thus $\ker \sigma = I_H$ for some saturated
hereditary subset $H \subseteq E^0$.  Let $\{ t_e, q_v \}$ be
the canonical Cuntz-Krieger $E$-family in $C^*(E)$.  Now
because
$\sigma(q_{v_0})$ is a rank 1 projection, it follows that
$q_{v_0} \notin \ker \sigma = I_H$ and thus $v_0 \notin H$. 
Since $H$ is hereditary this implies that for any $w \in H$ we
must have $w \ngeq v_0$.  Hence $H \subseteq H'$.

Now let $F := E / H$; that is, $F$ is the graph given by $F^0
:= E^0 \backslash H$ and $F^1 := \{ e \in E^1 : r(e) \notin H
\} $.  Then by \cite[Theorem 4.1]{BPRS} we see that $C^*(F)
\cong C^*(E) / I_H = C^*(E) / \ker \sigma$.  Thus we may
factor $\sigma$ as $\overline{\sigma} \circ p$ to get the
commutative diagram
\begin{equation}
  \notag
  \xymatrix{ \K \ar@{=}[d] \ar[r]^{i_E} &
C^*(E) \ar[d]^{\sigma} \ar[r]^{p} & C^*(F)
\ar[dl]^{\overline{\sigma}} \\  \K \ar[r]^i &
\B & }
\end{equation}
where $p$ is the standard projection and $\overline{\sigma}$
is the monomorphism induced by $\sigma$.  From the
commutativity of this diagram it follows that $p \circ i_E :
\K \rightarrow C^*(F)$ is injective.  Let $\{ s_e, p_v \}$ be
the canonical Cuntz-Krieger $F$-family in $C^*(F)$.  Also let
$I_{v_0}$ be the ideal in $C^*(E)$ generated by $q_{v_0}$,
and let $J_{v_0}$ be the ideal in $C^*(F)$ generated by
$p_{v_0}$. Using \cite[Corollary 2.2]{KPR} and the fact
that any path in $E$ with range $v_0$ is also a path in $F$,
we have that
\begin{align}
p(i_E(\K)) = & \ p(I_{v_0}) \notag \\
= & \ p ( \overline{\text{span}}\{ t_\alpha t_\beta^* : \alpha,
\beta \in E^* \text{ and } r(\alpha)=r(\beta) = v_0 \} ) \notag
\\ 
= & \ \overline{\text{span}}\{ p(t_\alpha t_\beta^*) : \alpha,
\beta \in E^* \text{ and } r(\alpha)=r(\beta) = v_0 \} \notag
\\ 
= & \ \overline{\text{span}}\{ s_\alpha s_\beta^* : \alpha,
\beta \in F^* \text{ and } r(\alpha)=r(\beta) = v_0 \} \notag
\\ 
= & \ J_{v_0}. \notag 
\end{align}
From the commutativity of the above diagram it follows
that $\overline{\sigma}$ is the (unique) map which makes the
diagram
\begin{equation}
  \notag
  \xymatrix{\K \ar[r]^{p \circ i_E}
\ar[dr]_{i}&
C^*(F) \ar[d]^{\overline{\sigma}} \\
& \B} 
\end{equation}
commute.  Since $\overline{\sigma}$ is injective, it follows
from \cite[Proposition 2.2.14]{WO} that $p(i_E(\K))=J_{v_0}$ is
an essential ideal in $C^*(F)$.

Now suppose that there exists $w \in F^0$ with $w \ngeq v_0$
in $F$.  Then for every $\alpha \in F^*$ with
$r(\alpha)=v_0$ we must have that $s(\alpha) \neq w$.  Hence
$p_w s_\alpha = 0$.  Since $J_{v_0} = \overline{\text{span}}
\{ s_\alpha s_\beta : \alpha, \beta \in F^* \text{ and }
r(\alpha)=r(\beta)=v_0 \}$ it follows that $p_w J_{v_0} = 0$. 
Since $p_w \neq 0$ this would imply that $J_{v_0}$ is not an
essential ideal.  Hence we must have that $w \ngeq v_0$ for
all $w \in F^0$.

Furthermore, if $\alpha \in F^*$ is a path with $s(\alpha)=w$
and $r(\alpha) = v_0$, then $\alpha \in E^*$.  So
if $w \ngeq v_0$ in $E$, then we must have that $w
\ngeq v_0$ in $F$.  Consequently, if $w \in H'$, then $w
\ngeq v_0$ in $E$, and we cannot have $w \in F^0$ because there is
a path in $F$ from every element of $F^0$ to $v_0$, and hence
a path in $E$ from every element of $F^0$ to $v_0$.  Thus $w
\notin F^0 := E^0 \backslash H$, and $w \in H$.  Hence
$H' \subseteq H$. \end{proof}

\begin{lemma}  
\label{wedgeequalsunion}
Let $G$ and $(E,v_0)$ be as before.  If $H_E$
is the inessential part of $E$, $H' := \{ v \in E^0 : v \ngeq
v_0 \}$, and $H_{v_0} := E^0 \backslash G^0$; then in $E$ we
have that
$$ H' \vee H_{v_0} = H_E \cup H_{v_0}.$$
\end{lemma}
\begin{proof}  We shall first show that $H_E \cup H_{v_0}$ is
a saturated hereditary subset of $E^0$.  To see that it is
hereditary, let $v \in H_E \cup H_{v_0}$.  If $e \in E^1$ with
$s(e)=v$, then one of two things must occur.  If $e \in G^1$,
then $s(e)=v$ must be in
$G^0$ and hence $v \in H_E$.  Since we know from Lemma
\ref{HEsaturatedhereditary} that
$H_E$ is a saturated hereditary subset of $G$, it follows that
$r(e) \in H_E \subseteq H_E \cup H_{v_0}$.  On the other hand,
if $e \notin G^1$, then $r(e) \notin G^0$, and hence $r(e) \in
H_{v_0} \subseteq H_E \cup H_{v_0}$.  Thus $H_E \cup H_{v_0}$
is hereditary.

To see that $H_E \cup H_{v_0}$ is saturated, let $v \notin H_E
\cup H_{v_0}$.  Then $v \in \overline{v_0}$ and $v \in \gamma$ for some
$\gamma \in \chi_G$ with the property that $\gamma \geq v_0$. 
Since maximal tails contain no sinks, there exists $e \in
G^1$ with $s(e) =v$ and $r(e) \in \gamma$.  But then $r(e)
\in \overline{v_0}$ and $r(e) \notin H_E$.  Since $e \in G^1$
this implies that $r(e) \notin H_E \cup H_{v_0}$.  Thus $H_E
\cup H_{v_0}$ is saturated.

Now since $H' \subset H_E$ we see that $H_E \cup H_{v_0}$ is a
saturated hereditary subset which contains $H' \cup H_{v_0}$. 
Thus $H' \vee H_{v_0} \subseteq H_E \cup H_{v_0}$.

Conversely, suppose that $v \in H_E \cup H_{v_0}$.  If $S$ is
any saturated hereditary subset of $E$ which contains $H' \cup
H_{v_0}$, then for every vertex $w \notin S$ we know that $w$
cannot be a sink, because if it were we would have $w \ngeq
v_0$.  Thus we may find an edge $e \in G^1$ with $s(e) = w$
and $r(e) \notin S$.  Furthermore, since $H' \cup H_{v_0}
\subseteq S$ we must also have that $r(e) \geq v_0$.  Thus if
$v \notin S$, we may produce an infinite path $\alpha$ in
$G$ with $s(\alpha) = v$ and
$s(\alpha_i) \geq v_0$ for all $i \in \N$.  If we let $\gamma
:= \{ w \in G^0 : w \geq s(\alpha_i) \text{ for some } i \in
\N \}$, then $\gamma \in \chi_G$ and $\gamma \geq v_0$.  Hence
$v \in \overline{v_0}$ and $v \notin H_E \cup H_{v_0}$ which is a
contradiction.  Thus we must have $v \in S$ for all saturated
hereditary subsets $S$ containing $H' \cup H_{v_0}$.  Hence $v
\in H' \vee H_{v_0}$ and $H_E \cup H_{v_0} \subseteq H' \vee
H_{v_0}$.  \end{proof}

\begin{lemma}
\label{kerneloftaudoesnotgotosink}
Let $G$ be a row-finite graph which satisfies Condition~(K). 
Also let $(E,v_0)$ be a 1-sink extension of
$G$.  If $\tau$ is the extension associated to
$E$, then $$\ker \tau = I_{H_E}.$$
\end{lemma}
\begin{proof}  Consider the following two cases.

\noindent \textsc{Case 1:}  The set $E^*(v_0)$ contains
finitely many elements.  

Then from the definition of the extension associated to $E$,
we have that $\tau = 0$.  However, if
$E^*(v_0)$ has only finitely many elements then $\gamma \ngeq
v_0$ for all $\gamma \in \chi_G$.  Hence $H_E = G^0$ and
$I_{H_E} = C^*(G)$.

\noindent \textsc{Case 2:}  The set $E^*(v_0)$ contains
infinitely many elements.  

Then $I_{v_0} \cong \K$, and from Lemma
\ref{kerneloftauisthewedgeofideals} we have that $\ker \tau =
\pi_E ( I_{v_0} \vee \ker \sigma)$.  Also Lemma
\ref{kernelofsigmaisvdontgotosink} implies that
$\ker \sigma = I_{H'}$.  Since $I_{v_0} = I_{H_{v_0}}$, we
see that from Lemma~\ref{wedgeequalsunion} that $I_{v_0} \vee
\ker \sigma = I_{H_{v_0}} \vee I_{H'} = I_{H_{v_0} \vee H'} =
I_{H_E \cup H_{v_0}}$.

Now if we let $\{s_e, p_v \}$ be the canonical Cuntz-Krieger
$G$-family in $C^*(G)$ and $\{ t_e, q_v \}$ be the canonical
Cuntz-Krieger $E$-family in $C^*(E)$, then 
$$ \ker \tau =  \pi_E(I_{H_E \cup H_{v_0}}) = \pi_E( \langle
\{ q_v : v \in H_E \cup H_{v_0} \} \rangle ) = \langle \{ p_v
: v \in H_E \} \rangle = I_{H_E}. $$
\end{proof}

\begin{lemma}
\label{pathstosinkarefinite}
Let $G$ be a row-finite graph which satisfies
Condition~(K), and let $(E,v_0)$ be a 1-sink extension of
$G$.  If $w \in H_E$, then 
\begin{equation*} \# \{ \alpha \in E^* : s(\alpha) = w \text{
and } r(\alpha)=v_0 \} < \infty. \end{equation*}
\end{lemma}
\begin{proof}  Suppose that there were infinitely many such
paths.  Then since $G$ is row-finite there must exist an edge
$e_1 \in G^1$ with $s(e_1)=w$ and with the property that
there exist infinitely many $\alpha \in E^*$ for which
$s(\alpha) = r(e_1)$ and $r(\alpha)=v_0$.  Likewise, there
exists an edge
$e_2 \in G^1$ with $s(e_2) = r(e_1)$ and with the property
that there are infinitely many $\alpha \in E^*$ for which
$s(\alpha)=r(e_2)$ and $r(\alpha)=v_0$.  Continuing in this
fashion we produce an infinite path $e_1e_2e_3\ldots$
with the property that $r(e_i) \geq v_0$ for all $i \in \N$. 
If we let $\gamma := \{ v \in G^0 : v \geq s(e_i) \text{ for
some }i \in \N \}$, then $\gamma \in \chi_G$ and $\gamma \geq
v_0$.  Since $w \in \gamma$, it follows that $w \in \overline{v_0}$ and
$w \notin H_E := E^0 \backslash \overline{v_0}$, which is a
contradiction. \end{proof}

\begin{definition}
Let $G$ be a row-finite graph which satisfies Condition~(K),
and let $(E,v_0)$ be a 1-sink extension of $G$.  Then $n_E
\in \prod_{H_E} \Z$ is the vector whose entries are given
by $$n_E(v) = \# \{ \alpha \in E^* : s(\alpha) = v \text{ and
} r(\alpha)=v_0 \} \hspace{.3in} \text{for $v \in H_E$}. $$
Note that the previous Lemma shows that $n_E(v) < \infty$ for
all $v \in H_E$.
\end{definition}

\begin{lemma}  Let $G$ be a row-finite graph which satisfies Condition~(K),
and let $(E,v_0)$ be a 1-sink extension of $G$.  If $v \in
H_E$ and $n_E(v) > 0$, then $A_G(v,v)=0$; that is, there does
not exist an edge $e \in G^1$ with $s(e)=r(e)=v$.
\label{noloopsatexits}
\end{lemma}
\begin{proof} If there was such an edge $e \in G^1$, then
$\gamma = \{ w \in G^0 : w \geq v \}$ would be a maximal tail
and since $n_E(v) > 0$ it would follow that $\gamma \geq
v_0$.  Since $v \in \gamma$ this implies that $v \in
\overline{v_0}$ which contradicts the fact that $v \in H_E :=
G^0 \backslash \overline{v_0}$.
\end{proof}

\begin{lemma} 
\label{rankofedgewithrangeinhe}
Let $G$ be a row-finite graph which satisfies Condition~(K),
and let $(E,v_0)$ be a 1-sink extension of $G$.  Also let $\{
t_e, q_v \}$ be the canonical Cuntz-Krieger $E$-family in
$C^*(E)$.  If $e \in G^1$ and $r(e) \in H_E$, then 
$$\rank \sigma (t_e) = n_E(r(e)).$$
\end{lemma}
\begin{proof}  If $n_E(r(e)) = 0$, then $r(e) \ngeq v_0$ and
by Lemma \ref{kernelofsigmaisvdontgotosink} we have $\sigma (
q_{r(e)})=0$.  Since $\sigma (t_e)$ is a partial isometry
$\rank \sigma(t_e) = \rank \sigma (t_e^*t_e) = \rank \sigma
(q_{r(e)}) = 0$.  Therefore we need only consider the case
when $n_E(r(e)) > 0$.

Let $B_E^1$ denote the boundary edges of $E$.  Also let $k_e
:= \text{max} \{ | \alpha | : \alpha
\in E^*, s(\alpha)=r(e), \text{ and } r(\alpha) \in B_E^1
\}$.  By Lemma
\ref{pathstosinkarefinite} we see that $k_e$ is finite.  We
shall prove the claim by induction on $k_e$.

\noindent Base Case: $k_e = 0$.  Let $e_1, e_2, \ldots , e_n$
be the boundary edges of $E$ which have source $r(e)$.  Then
it follows from \cite[Lemma 6.9]{Tom} that for $1 \leq i
\leq n$ $$\rank \sigma(t_{e_i}) = \# Z(r(e_i),v_0)$$ where
$Z(r(e_i),v_0)$ is the set of paths from $r(e_i)$ to $v_0$. 
Also if $f
\in G^1$ is an edge with $s(f) = r(e)$, then because
$n_E(r(e)) > 0$ Lemma \ref{noloopsatexits} implies that $r(f)
\neq r(e)$.  Furthermore, since $k_e =0$ we must have that
$r(f) \ngeq v_0$.  Therefore, just as before we must have
$\rank \sigma (t_f) = 0$.  Now since the projections $\{
t_ft_f^* : f \in E^1 \text{ and } s(f)=r(e)
\}$ are mutually orthogonal, we see that 
\begin{align}
\rank \sigma (t_e) = & \ \rank \sigma (t_e^*t_e) \notag \\
= & \ \rank \sum_{ f \in E^1 \atop s(f) = r(e) } \sigma
(t_ft_f^*) \notag \\
= & \ \rank \sigma (t_{e_1}) + \ldots + \rank \sigma (t_{e_n})
+ \sum_{ f \in G^1 \atop s(f) = r(e) } \rank \sigma
(t_ft_f^*) \notag \\
= & \ \# Z(r(e_1),v_0) + \ldots + \# Z(r(e_n),v_0). \notag \\
= & \ n_E(r(e)). \notag
\end{align}

\noindent Inductive Step:  Assume that the claim holds for all
edges $f$ with $k_f \leq m$.  We shall now show that the
claim holds for edges $e \in G^1$ with $k_e = m+1$.  Let $e_1
, e_2, \ldots, e_n$ be the exits of $E$ with source $r(e)$. 
As above we have that $\rank \sigma (t_{e_i}) = \#
Z(r(e_i),v_0)$ for all $1 \leq i \leq n$.  Now if $f \in G^1$
is any edge with $s(f) = r(e)$, then Lemma
\ref{noloopsatexits} implies that $r(f) \neq r(e)$.  Thus we
must have that $k_f \leq k_e - 1$, and by the induction
hypothesis $\rank \sigma (t_f) = n_E(r(f))$.  Furthermore,
since the projections $\{ t_ft_f^* : f \in E^1
\text{ and } s(f)=r(e) \}$ are mutually orthogonal, we see
that 
\begin{align}
\rank \sigma (t_e) = & \ \rank \sigma (t_e^*t_e) \notag \\
= & \ \rank \sum_{ f \in E^1 \atop s(f) = r(e) } \sigma
(t_ft_f^*) \notag \\
= & \ \rank \sigma (t_{e_1}) + \ldots + \rank \sigma (t_{e_n})
+ \sum_{ f \in G^1 \atop s(f) = r(e) } \rank \sigma
(t_ft_f^*) \notag \\
= & \ \# Z(r(e_1),v_0) + \ldots + \# Z(r(e_n),v_0) + \sum_{ f
\in G^1 \atop s(f) = r(e) } n_E(r(f))
\notag \\
= & \ n_E(r(e)). \notag
\qed
\end{align}
\renewcommand{\qed}{}
\end{proof}

Let $G$ be a row-finite graph which satisfies Condition~(K)
and let $(E,v_0)$ be a 1-sink extension of $G$.  If $H_E :=
G^0 \backslash \overline{v_0}$ is the inessential part of
$E$, then since $H_E$ is a saturated hereditary subset of $G$
we may form the graph $F := G / H_E$ given by $F^0 :=
G^0 \backslash H_E$ and $F^1 := \{ e
\in G^1 : r(e) \notin H_E \}$.  With respect to the
decomposition $G^0 = \overline{v_0} \cup H_E$ the vertex
matrix $A_G$ of $G$ will then have the form
$$A_G = \begin{pmatrix} A_F & X \\ 0 & C \end{pmatrix}$$
where $A_F$ is the vertex matrix of the graph $F$. 

Furthermore, if $\tau : C^*(G) \rightarrow \Q$ is the Busby
invariant of the extension associated to $E$, then by Lemma
\ref{kerneloftaudoesnotgotosink} we know that $\ker \tau =
I_{H_E}$.  Hence $C^*(G)/ \ker \tau \cong C^*(F)$ by
\cite[Theorem 4.1]{BPRS} and we may factor $\tau$ as
$\overline{\tau}
\circ p$
\begin{equation}
  \notag
  \xymatrix{C^*(G) \ar[r]^{p}
\ar[d]_{\tau}&
C^*(F) \ar[dl]^{\overline{\tau}} \\
\Q &  } 
\end{equation}
where $p$ is the standard projection and $\overline{\tau}$ is
the monomorphism induced by $\tau$.  Note that since
$\overline{\tau}$ is injective it is an essential
extension of $C^*(F)$.  Furthermore, with respect to the
decomposition $G^0 = \overline{v_0} \cup H_E$ the $\W$ vector of $E$ will
have the form $\wvec_E = \left( \begin{smallmatrix} \wvec_E^1
\\ \wvec_E^2 \end{smallmatrix} \right)$.

\begin{lemma}
If $d : \ext (C^*(F)) \rightarrow \coker (B_F -I)$ is the
Cuntz-Krieger map, then 
$$d(\overline{\tau}) = [x]$$
where $[x]$ denotes the class in $\coker (B_F-I)$ of the vector
$x \in \prod_{F^1} \Z$ given by $x(e) := \wvec_E^1(r(e)) +
(Xn_E)(r(e))$ for all $e \in F^1$. 
\label{doftaubarisxplusjunk}
\end{lemma}
\begin{proof}  Notice that because of the way $H_E$ was
defined, $F$ will have no sinks.  Also note that the diagram 
\begin{equation}
  \notag
  \xymatrix{ \K \ar@{=}[d] \ar[r]^{i_E} &
C^*(E) \ar[d]^{\sigma} \ar[r]^{\pi_E} & C^*(G) \ar[r]^p
\ar[d]^\tau & C^*(F)
\ar[dl]^{\overline{\tau}} \\  \K \ar[r]^i &
\B \ar[r]^\pi & \Q & }
\end{equation}
commutes.  Let $\{t_e, q_v \}$ be the
canonical Cuntz-Krieger $E$-family in $C^*(E)$. For each $e
\in F^1$ let $$H_e := \im \sigma (t_et_e^*)$$ and for each
$v \in F^0$ let $$H_v := \bigoplus_{ e \in F^1 \atop
s(e)=v } H_e.$$  Also for each $v \in F^0$ define $P_v$ to
be the projection onto $H_v$ and for each $e \in F^1$ define
$S_e$ to be the partial isometry with initial space $H_{r(e)}$
and final space $H_e$.  Then $\{ S_e, P_v \}$ is a
Cuntz-Krieger $F$-family in $\B$.  If we let $\{s_e, p_v \}$
be the canonical Cuntz-Krieger $F$-family in $C^*(F)$, then by
the universal property of $C^*(F)$ there exists a
homomorphism $\tilde{t} : C^*(F) \rightarrow \B$ such that
$\tilde{t} (s_e) = S_e$ and $\tilde{t}(p_v)=P_v$.  Let
$t := \pi \circ \tilde{t}$.  Since $G$ satisfies
Condition~(K), it follows that the quotient $F := G / H_E$
also satisfies Condition~(K).  Because $t (p_v) \neq 0$ for
all $v \in F$ this implies that $\ker t = 0$ and $t$ is an
essential extension of $C^*(F)$.

Because $\sigma (t_e)$ is a lift of $\overline{\tau} (s_e)$
for all $e \in F^1$ we see that $\ind_{E_e} \overline{\tau}
(s_e) t (s_e)$ equals the Fredholm index of $\sigma (t_e)
S_e^*$ in $H_e$.  Since $S_e^*$ is a partial isometry with
initial space $H_e$ and final space $H_{r(e)} \subseteq \im
\sigma (q_{r(e)})$, and since $\sigma (t_e)$ is a partial
isometry with initial space $\im \sigma (q_{r(e)})$ and final
space $H_e$, it follows that $$\text{dim} (\ker (\sigma (t_e)
S_e^*)) = 0.$$
Also, $\sigma (t_e^*)$ is a partial isometry with
initial space $H_e$ and final space $\im \sigma (q_{r(e)})$,
and $S_e$ is a partial isometry with initial space
$H_{r(e)}$ and final space $H_e$.  Because $q_{r(e)} = \sum_{
\{ f \in E^1 : s(f) = r(e) \} } t_ft_f^*$ we see that $$\im
\sigma (q_{r(e)}) = H_{r(e)} \oplus \bigoplus_{ f \in E^1
\backslash F^1 \atop s(f)=r(e) } \im \sigma (t_f t_f^*) $$
Thus $$\text{dim} (\ker ( S_e \sigma (t_e^*))) = \sum_{ f
\in E^1 \backslash F^1 \atop s(f) = r(e) } \rank \sigma
(t_ft_f^*) = \sum_{ f \in E^1 \backslash F^1 \atop s(f)=r(e)
} \rank \sigma (t_f).$$
Now if $f$ is any boundary edge of $E$, then by
\cite[Lemma 6.9]{Tom} we have that $\rank
\sigma (t_f) = \# Z((r(f),v_0)$ where $Z((r(f),v_0)$ is the
set of paths in $F$ from $r(f)$ to $v_0$.  Also if
$f \in G^1$ is an edge with $r(f) \in H_E$, then $\rank \sigma
(t_f) = n_E(r(f))$ by Lemma \ref{rankofedgewithrangeinhe}. 
Therefore,
\begin{align}
\text{dim} (\ker (S_e \sigma (t_e^*))) = & \ \sum_{ \text{$f$
is a boundary edge} \atop s(f) = r(e) } \rank \sigma (t_f) +
\sum_{ r(f) \in H_E \atop s(f) = r(e)} \rank \sigma (t_f)
\notag \\ = & \ \sum_{ \text{$f$
is a boundary edge} \atop s(f) = r(e) } \# Z(r(f),v_0) +
\sum_{r(f) \in H_E
\atop s(f) = r(e)} n_E(r(f)) \notag \\
= & \ \wvec_E(r(e)) + \sum_{w \in H_E} X(r(e),w) n_E(w). \notag
\end{align}
Thus \begin{align} d_{\overline{\tau},t} (e) & = \ -
\ind_{E_e} \overline{\tau} (s_e) t(s_e^*) \notag \\ 
& = \ \wvec_E(r(e)) + \sum_{w \in H_E} X(r(e),w) n_E(w)
\notag \\
& =  \ \wvec_E^1(r(e)) +(Xn_E)(r(e)). \notag 
\qed
\end{align}
\renewcommand{\qed}{}
\end{proof}

\begin{lemma}
If $\wmap : \ext (C^*(F)) \rightarrow \coker (A_F-I)$ is the
$\W$ map, then $$\wmap (\overline{\tau}) = [ \wvec_E^1 + Xn_E
]$$ where $[\wvec_E^1 +Xn_E]$ denotes the class of the vector
$\wvec_E^1 + X n_E$ in $\coker (A_F-I)$.
\label{valueofwmapiswvecplusjunk}
\end{lemma}
\begin{proof}  By definition $\wmap := \overline{S_F} \circ
d$.  From Lemma \ref{doftaubarisxplusjunk} we see that
$d(\overline{\tau}) = [x]$, where $x(e) = \wvec_E^1(r(e)) +
(Xn_E) (r(e))$ for all $e \in F^1$.  Therefore, $\wmap (
\overline{\tau} )$ is equal to the class $[y]$ in
$\coker(A_F-I)$ where $y \in \prod_{F^0} \Z$ is the vector
given by $y := S_F (x)$.  Hence for all $v
\in F^0$ we have that  $$y (v) = (S_F (x))(v) = \sum_{ e \in
F^1 \atop s(e) = v } x(e) = \sum_{ e \in F^1 \atop s(e) = v}
\wvec_E^1 (r(e)) + (Xn_E) (r(e))$$ 
and thus for all $v \in F^0$ we have that
\begin{align}
& \ y(v) - \left( \wvec_E^1(v) +
(Xn_E)(v)
\right) \notag \\
= & \ \left( \sum_{ e \in F^1 \atop s(e) = v}
\wvec_E^1(r(e)) + (Xn_E)(r(e)) \right) - \left( \wvec_E^1(v) +
(Xn_E)(v) \right) \notag \\
= & \ \left( \sum_{w \in F^0} A_F(v,w) \left( \wvec_E^1(w)
+(Xn_E)(w) \right) \right) - \left( \wvec_E^1(v) + (Xn_E)(v)
\right). \notag
\end{align}
Hence $y - (\wvec_E^1 + Xn_E) = (A_F-I)
(\wvec_E^1 +Xn_E)$, and
$\wmap (\tau) = [y]=[\wvec_E^1+Xn_E]$ in $\coker
(A_F-I)$. \end{proof}

\begin{remark}Let $G$ be a row-finite graph which
satisfies Condition~(K) and let $(E_1,v_1)$ and $(E_2,v_2)$
be 1-sink extensions of $G$.  If $\overline{v_1} =
\overline{v_2}$, then we may let $H := H_{E_1} = H_{E_2}$ and
form the graph $F := G / H$ given by $F^0 := G^0 \backslash
H$ and $F^1 := \{ e \in G^1 : r(e) \notin H \}$).  Then with
respect to the decomposition $G^0 = (G^0 \backslash H) \cup
H$, the vertex matrix of $G$ has the form $$A_G =
\begin{pmatrix} A_F & X \\ 0 & C \end{pmatrix}$$ where $A_F$
is the vertex matrix of
$F$.  Also with respect to this decomposition, the $\W$
vectors of $E_1$ and $E_2$ have the form $\wvec_{E_1} = \left(
\begin{smallmatrix} \wvec_{E_1}^1 \\ \wvec_{E_1}^2
\end{smallmatrix} \right)$ and $\wvec_{E_2} = \left(
\begin{smallmatrix} \wvec_{E_2}^1 \\ \wvec_{E_2}^2
\end{smallmatrix} \right)$. \label{decompremark} \end{remark}

For $i \in \{ 1, 2 \}$, let $n_{E_i} \in
\prod_H \Z$ denote the vector given by $n_{E_i}(v) = \# \{
\alpha \in E_i^* : s(\alpha) = v \text{ and } r(\alpha) = v_i
\}$. 

\begin{theorem}  Let $G$ be a row-finite graph which
satisfies Condition~(K), and let $(E_1,v_1)$ and $(E_2,v_2)$
be 1-sink extensions of $G$.  Using the notation in
Remark~\ref{decompremark}, we have that one of the
$C^*(E_i)$'s may be $C^*(G)$-embedded into the other if and
only~if \begin{enumerate} \item $\overline{v_1} =
\overline{v_2}$
\item $[\wvec_{E_1}^1 + X n_{E_1}] = [\wvec_{E_2}^1 + X
n_{E_2}] \text{ in $\coker(A_F-I)$}.$
\end{enumerate} 
\label{charctGequivforgeneralexten}
\end{theorem}
\begin{proof}  It follows from Lemma \ref{primconverse} that
if one of the $C^*(E_i)$'s is $C^*(G)$-embeddable in the
other, then $\overline{v_1} =
\overline{v_2}$.  Thus we may let $H := H_{E_1} = H_{E_2}$
and form the graph $F := G / H$ as discussed in Remark
\ref{decompremark}.  

If we let $\tau_1$ and $\tau_2$ be the Busby invariants of the
extensions associated to $E_1$ and $E_2$, then it follows from
Lemma \ref{kerneloftaudoesnotgotosink} that $\ker
\tau_1 = \ker \tau_2 = I_H$.  Thus for each $i \in \{ 1, 2
\}$, we may factor $\tau_i$ as $\tau_i = \overline{\tau}_i
\circ p$
\begin{equation}
  \notag
  \xymatrix{C^*(G) \ar[r]^{p}
\ar[d]_{\tau_i}&
C^*(F) \ar[dl]^{\overline{\tau}_i} \\
\Q &  } 
\end{equation}
where $p$ is the standard projection and $\overline{\tau}_i$
is the monomorphism induced by $\tau_i$.

It then follows from Theorem \ref{GequiviffCKequiv} that
one of the $C^*(E_i)$'s may be $C^*(G)$-embedded into the
other if and only if $\tau_1$ and
$\tau_2$ are CK-equivalent.  Since $\tau_i = \overline{\tau}_i
\circ p$ we see that $\tau_1$ and $\tau_2$ are CK-equivalent if
and only if $\overline{\tau}_1$ and $\overline{\tau}_2$ are
CK-equivalent.  Furthermore, since $\overline{\tau}_1$ and
$\overline{\tau}_2$ are essential extensions we see from
Corollary \ref{Gequiviffwsequiv} that $\overline{\tau}_1$ and
$\overline{\tau}_2$ are CK-equivalent if and only if
$\overline{\tau}_1$ and $\overline{\tau}_2$ are equal in $\ext
(C^*(F))$.  If $\wmap : \ext (C^*(F)) \rightarrow \coker
(A_F-I)$ is the $\W$ map, then this will occur if and only if
$\wmap (\overline{\tau}_1) = \wmap ( \overline{\tau}_2)$, and
by Lemma \ref{valueofwmapiswvecplusjunk} we see that this
happens if and only if $[\wvec_{E_1}^1 + Xn_{E_1}] =
[\wvec_{E_2}^1 + Xn_{E_2}]$ in $\coker (A_F-I)$. \end{proof}

\begin{remark}  Note that when $E_1$ and $E_2$ are
both essential we have $\overline{v_1} = \overline{v_2} =
G^0$ and $H= \emptyset$.  In this case
$F =G$, $X$ is empty, and
$\wvec_{E_i}^1 = \wvec_{E_i}$ for $i =1,2$.  Thus the result
for essential extensions in Theorem
\ref{charactofGequivforessential} is a special case of the
above theorem.   

In addition, we see that the above theorem gives a
method of determining $C^*(G)$-embeddability from
basic calculations with data that can be easily read of from
the graphs.  To begin, the condition that $\overline{v_1}
= \overline{v_2}$ can be checked simply by looking at $E_1$
and $E_2$.  In addition, the set $H$, the matrices $A_F$ and
$X$, and the vectors $\wvec_{E_i}^1$ and $n_{E_i}$ for
$i=1,2$ can easily be read off from the graphs $G$, $E_1$,
and $E_2$.  Finally, determining whether $[\wvec_{E_1}^1 + X
n_{E_1}] = [\wvec_{E_2}^1 + X n_{E_2} ]$ in $\coker(A_F-I)$
amounts to ascertaining whether $(\wvec_{E_1}^1 -
\wvec_{E_2}^1 ) + ( X(n_{E_1}-n_{E_2}) ) \in \im (A_F-I)$, a
task which reduces to checking whether a system of linear
equations has a solution.
\end{remark}

We now mention an interesting consequence of the above theorem.
\begin{definition}
Let $G$ be a row-finite graph which satisfies Condition~(K),
and let $(E,v_0)$ be a 1-sink extension of $G$.  We say that
$E$ is \emph{totally inessential} if
$\overline{v_0}=\emptyset$; that is, if $\{ \gamma \in \chi_G
: \gamma \geq v_0 \} = \emptyset$.
\end{definition}

\begin{corollary}
Let $G$ be a row-finite graph which satisfies Condition~(K).
If $(E_1,v_1)$ and $(E_2,v_2)$ are 1-sink extensions of
$G$ which are totally inessential, then one of the
$C^*(E_i)$'s may be $C^*(G)$-embedded into the other.
\label{totallyinessentialimliesGequiv}
\end{corollary}

\begin{proof}  Using the notation established in Remark
\ref{decompremark} and the proof of Theorem
\ref{charctGequivforgeneralexten}, we see that if
$E_1$ and $E_2$ are totally inessential, then $H = G^0$. 
Hence $F = \emptyset$ and $\overline{\tau}_1 =
\overline{\tau}_2=0$.  Thus $\overline{\tau}_1$ and
$\overline{\tau}_2$ are trivially CK-equivalent.  Hence
$\tau_1$ and $\tau_2$ are CK-equivalent and it follows from
Theorem \ref{GequiviffCKequiv} that one of the $C^*(E_i)$'s
can be $C^*(G)$-embedded into the other.

Alternatively, we see that if $E_1$ and $E_2$ are totally
inessential, then $F = \emptyset$, and provided that we
interpret the condition that $[\wvec_{E_1}^1 + X n_{E_1}] =
[\wvec_{E_2}^1 + X n_{E_2}]$ in
$\coker(A_F-I)$ as being vacuously satisfied, the
previous theorem implies that one of the $C^*(E_i)$'s
can be $C^*(G)$-embedded into the other.  \end{proof}

\begin{remark}
The case when $E_1$ and $E_2$ are both essential and the case
when $E_1$ and $E_2$ are both inessential can be thought of as
the degenerate cases of Theorem
\ref{charctGequivforgeneralexten}.  The first occurs when
$\overline{v_0} = G^0$ and $H = \emptyset$, and the second
occurs when $\overline{v_0} = \emptyset$ and $H = G^0$.
\end{remark}

\begin{example}  Let $G$ be the graph
\begin{equation}
  \notag
  \xymatrix{\cdots \ar[r]^{e_{-2}} & v_{-1} \ar[r]^{e_{-1}} &
v_0 \ar[r]^{e_0} & v_1 
    \ar[r]^{e_1}  & v_2 \ar[r]^{e_2} & v_3 \ar[r]^{e_3} &
\cdots} 
\end{equation}  
Note that $C^*(G) \cong \K$.  Since $G$ has precisely one
maximal tail $\gamma := G^0$, we see that if $E$ is any 1-sink extension of $G$, then $E$ will either be essential or
totally inessential.  Furthermore, one can check that $A_G-I
: \prod_{G^0} \Z \rightarrow \prod_{G^0} \Z$ is surjective. 
Thus if $E_1$ and $E_2$ are two essential 1-sink extensions
of $G$, we will always have that $[\wvec_{E_1}] =
[\wvec_{E_2}]$ in $\coker(A_G-I)$.  In light of Theorem
\ref{charactofGequivforessential} and Corollary
\ref{totallyinessentialimliesGequiv} we see that if $E_1$ and
$E_2$ are two 1-sink extensions of $G$, then one of the
$C^*(E_i)$'s can be $C^*(G)$-embedded in to the other if and
only if they are both essential or both totally inessential. 
\end{example}

We end with an interesting observation.  Note
that the statement of the result in Theorem
\ref{charctGequivforgeneralexten} involves the
$\wvec_{E_i}^1$ terms from the $\W$ vectors, but does not make
use of the $\wvec_{E_i}^2$ terms.  If for each
$i \in \{ 1, 2 \}$ we let $B^0_{E_i}$
denote the boundary vertices of $E_i$, then we see that the
nonzero terms of $\wvec_{E_i}^2$ are those entries which
correspond to the elements of $B^0_{E_i}
\cap H$.  Furthermore, if $v \in B^0_{E_i} \cap H$ and there
is a path from $F^0$ to $v$, then the value of
$\wvec_{E_i}^2(v)$ will affect the value of $n_{E_i}$. 
However, if there is no path from $F^0$ to $v$, then the
value of $\wvec_{E_i}^2(v)$ will be irrelevant to the value
of $n_{E_i}$.

Therefore, whether one of the $C^*(E_i)$'s may be
$C^*(G)$-embedded onto the other depends on two things:
the number of boundary edges at vertices in
$F^0$ (which determine the value of the $\wvec_{E_i}^1$'s),
and the number of boundary edges at vertices in $H$ which can
be reached by $F^0$ (which determine the value of the
$n_{E_i}$'s).  The boundary edges whose sources are elements
of $H$ that cannot be reached by $F^0$ will not matter.

\end{document}